\documentclass[onefignum,onetabnum]{siamart171218}

\usepackage{amsmath}
\usepackage{lineno}
\usepackage{amsfonts}
\usepackage{algorithm}
\usepackage{algorithmicx}
\usepackage{empheq}
\usepackage{amsfonts}
\usepackage{subfig}
\usepackage{enumitem}
\usepackage{mathtools,cases}
\usepackage{mathtools,mathptmx}
\usepackage{mathtools}
\usepackage{amssymb}
\usepackage{appendix}

\newsiamthm{remark}{Remark}
\headers{Convexification Method for a 1D Wave Equation}{A. Smirnov, M. Klibanov, A. Sullivan and L. Nguyen}

\begin{document}

\title{Convexification for an Inverse Problem \\
for a 1D Wave Equation with Experimental Data\thanks{
Submitted to the editors DATE. 
\funding{The workof AVS and MVK was supported
by US Army Research Laboratory and US Army Research Office grant W911NF-19-1-0044.}}}

\author{Alexey V. Smirnov\thanks{Department of Mathematics and Statistics, University of North Carolina at Charlotte, Charlotte, NC, 28223 (\email{asmirno2@uncc.edu}, \email{mklibanv@uncc.edu} (corresponding author)),
} \and Michael V. Klibanov\footnotemark[2] \and Anders Sullivan\thanks{US Army Research Laboratory, 2800 Powder Mill Road Adelphy, MD 20783-1197 \newline (\email{anders.j.sullivan.civ@mail.mil}, \email{lam.h.nguyen2.civ@mail.mil}).   
\newline
} \and Lam Nguyen\footnotemark[3]}
\maketitle

\begin{abstract}
The forward problem here is the Cauchy problem for a 1D hyperbolic PDE with
a variable coefficient in the principal part of the operator. That
coefficient is the spatially distributed dielectric constant. The inverse
problem consists of the recovery of that dielectric constant from
backscattering boundary measurements. The data depend on one variable, which is time. To address this problem, a new version of the convexification
method is analytically developed. The theory guarantees the global
convergence of this method. Numerical testing is conducted for both
computationally simulated and experimental data. Experimental data, which
are collected in the field, mimic the problem of the recovery of the
spatially distributed dielectric constants of antipersonnel land mines and
improvised explosive devices.
\end{abstract}


\noindent AMS Classification 35R30 

\section{Introduction}

\label{sec:1}

In this paper, we first develop a new globally convergent numerical method for a Coefficient Inverse Problem (CIP) for a 1D hyperbolic equation. This is a version of the so-called \emph{convexification} method, which has been
actively developed by the second coauthor and his coauthors for the past several years, see, e.g. \cite{KlibanovNik:ra2017,khoa2019convexification,KlibanovTimonov:u2004,Klibanov:sjma1997,klibanov1995uniform,klibanov2019backscatter,klibanovzhang2020convexification,klibanov2019convexification,klibanov2018new,smirnov2020convexification}
for some samples of those publications. First, we test this method on
computationally simulated data. Next, we test it on experimental time
dependent data collected by the US Army Research Laboratory (ARL) using
their forward looking radar \cite{nguyen2007obstacle}. These data were
collected in the field which is a more difficult case than the collection in
a laboratory. A significant challenge here is that targets were surrounded
by clutter.

That radar was built to image and identify flash explosive-like targets,
mainly antipersonnel land mines and improvised explosive devices. Those
targets can be located both in air and under the ground. In the latter case,
the burial depth does not exceed 10 centimeters. Currently ground
penetrating radars rely only on the information about the energy of the
backscattering signals. In our case, however, we computationally estimate
dielectric constants of those targets. We hope that estimates of dielectric
constants of explosive-like targets might serve in the future classification
algorithms as an additional parameter to the currently used ones. The use of
an additional parameter, in turn might decrease the current false alarm rate.

Although this research group has a number of publications where these experimental data are treated, in each of them either Laplace or Fourier transform of the experimental data with respect to the time $t$ is considered \cite{karchevsky2013krein,KlibanovLoc:ipi2016,klibanov2018new,klibanov2017globally,Kuzhuget:IP2012}%
. But since our data are actually time resolved ones, then we work here
directly in the time domain.

Given a CIP, we call a numerical method for it \emph{globally convergent},
if the theory rigorously guarantees that one can obtain at least one point
in a sufficiently small neighborhood of the exact solution \emph{without} an
assumption that the starting point of iterations is located in this
neighborhood. The global convergence issue is obviously important since it
is very rare in applications when a good first approximation for the
solution is available \emph{a priori}. 

Any CIP is both nonlinear and ill-posed. These factors cause the well known
phenomenon of multiple local minima and ravines of conventional
Tikhonov-like least squares functionals, see, e.g. \cite{ScalesSmithFischerLjcp1992} for a convincing example of this phenomenon.
Conventional numerical methods for CIPs rely on the minimization of such
functionals, see, e.g. \cite{chavent2010nonlinear} for some examples with a
study of some experimental data in \cite{goncharsky2019method,goncharsky2019low}.

Unlike conventional approaches, the concept of the convexification is
designed with the goal to avoid local minima and ravines, see section 4.1
for a brief outline of the convexification. The first publications about
this concept were in 1995 and 1997 \cite{Klibanov:sjma1997,klibanov1995uniform}. However, numerical studies were not
conducted at that time since there were a number of theoretical issues which
needed to be clarified to pave the way for numerical testing, although, see
some numerical examples in the book \cite{KlibanovTimonov:u2004}. Those
issues were clarified more recently in the paper \cite{KlibanovNik:ra2017}.
This publication has generated a number of papers with numerical studies. As
some examples of those we mention \cite{khoa2019convexification,klibanov2019backscatter,klibanovzhang2020convexification,klibanov2019convexification,klibanov2018new,klibanov2017globally,smirnov2020convexification}.
In particular, in \cite{klibanov2019convexification} the convexification is applied to a CIP with single measurement data for a 3D hyperbolic PDE.

To work with the experimental data of ARL in the time domain, one needs to
have a numerical method for a CIP for a 1D hyperbolic PDE. However, the
convexification for CIPs for the 1D hyperbolic PDEs was not developed until
the very recent work \cite{smirnov2020convexification}. Even though the CIP
of the current paper is reduced to the same CIP as the one in \cite{smirnov2020convexification} for the equation $u_{tt}=u_{xx}+p\left(
x\right) u,$ the method of this paper is significantly different from the
one of \cite{smirnov2020convexification}. Indeed, in \cite{smirnov2020convexification} an integral differential equation with the
Volterra-like integrals is obtained under the condition that the unknown
coefficient $p\left( x\right) \geq 0.$ Unlike this, in the current paper we
replace those Volterra-like integrals with a non local condition. This
allows us to avoid imposing the non-negativity condition on $p\left(
x\right) .$ The latter, in turn enables us not to impose some quite
restrictive extra conditions on the target unknown coefficient $c\left(
y\right) ,$ see (\ref{2.3}) and (\ref{2.8}).

We now refer to some other numerically implemented globally convergent numerical methods for CIPs for hyperbolic PDEs. First, this is the paper of De Hoop, Kepley and Oksanen \cite{de2018recovery} for a 3D hyperbolic CIP with Dirichlet-to-Neumann map data. Next, this is a series of works of Kabanikhin with coauthors who have computationally implemented the Gelfand-Levitan method \cite{karchevsky2013krein,khoa2019convexification} for both 1D and 2D cases. We also mention here works of Baudouin, de~Buhan, Ervedoza and Osses \cite{baudouin2017convergent,baudouin2020carleman}, where a different version of the convexification for CIPs for hyperbolic PDEs in the $n-$D case is developed, also see more recent works of Boulakia, de Buhan and Schwindt \cite{boulakia2019numerical} and Le and Nguyen \cite{le2019convergent}, where the idea of \cite{baudouin2017convergent} is developed further to apply to some nonlinear inverse problems for parabolic PDEs.  The common property of these works and the above cited works on the convexification is that both substantially use Carleman estimates and the resulting numerical methods converge globally in both cases.

However, there is a significant difference between \cite{baudouin2017convergent,baudouin2020carleman} and the above cited
publications on the convexification \cite{khoa2019convexification,KlibanovTimonov:u2004,Klibanov:sjma1997,klibanov1995uniform,klibanov2019backscatter,klibanov2019convexification,klibanov2018new,klibanov2017globally,smirnov2020convexification}. More precisely, \cite{baudouin2017convergent,baudouin2020carleman} work
under an assumption of the 
Bukhgeim-Klibanov method \cite{BukhgeimKlibanov:smd1981}. This assumption is
that one of initial conditions is not vanishing in the entire domain of
interest. The original version of the convexification uses this assumption
in \cite{BeilinaKlibanovBook,klibanovzhang2020convexification}. On the other
hand, in publications \cite{khoa2019convexification,KlibanovTimonov:u2004,Klibanov:sjma1997,klibanov1995uniform,klibanov2019backscatter,klibanov2019convexification,klibanov2018new,klibanov2017globally,smirnov2020convexification}, so as in this paper, either the initial condition, or a corresponding
right hand side of a PDE is vanishing.

The convexification has significant roots in the idea of the paper \cite%
{BukhgeimKlibanov:smd1981} (1981), in which the tool of Carleman estimates
was introduced in the field of CIPs for the first time. The original goal of 
\cite{BukhgeimKlibanov:smd1981} was limited to proofs of global uniqueness
theorems for multidimensional CIPs. Therefore, the convexification can be
regarded as an extension of the concept of \cite{BukhgeimKlibanov:smd1981}
from the purely theoretical uniqueness topic to a more applied topic of
globally convergent numerical methods for CIPs. Since 1981, many works of
many authors have been devoted to a variety of applications of the method of 
\cite{BukhgeimKlibanov:smd1981} to proofs of uniqueness and stability
results for many CIPs. Since the current paper is not a survey of these
works, we now refer here only to the books \cite%
{BeilinaKlibanovBook,BellassouedYamamoto:SpKK2017,KlibanovTimonov:u2004} and
the survey \cite{Klibanov:jiipp2013}.

All functions below are real valued ones. In section 2, we formulate both
forward and inverse problems. In section 3, we derive a quasilinear 1D
hyperbolic PDE with a non local condition. In section 4, we introduce a
weighted Tikhonov-like functional, which is the main subject of the
convexification. In section 5, we formulate theorems about that functional,
ensuring the global convergence of the resulting gradient projection method.
Section 6 contains the proofs of theorems formulated in section 5. In
section 7, we describe the algorithms, as well as accompanying procedures
used to obtain the numerical results for both computationally simulated and
experimental data.

\section{Statements of Forward and Inverse Problems}

\label{sec:2} Let $\overline{c}>1$ be a number and the function $c\left(
y\right) \in C^{3}\left( \mathbb{R}\right) $ has the following properties%
\begin{align}
& c\in \left[ 1,\overline{c}\right] ,\quad \overline{c}=const>1,  \label{2.1}
\\
& c\left( y\right) =1,\quad y\in (-\infty ,0]\cup \lbrack 1,+\infty ).
\label{2.2}
\end{align}%
Physically $c\left( y\right) =n^{2}\left( y\right) $ is the spatially
distributed dielectric constant, where $n\left( y\right) $ is the refractive
index. In acoustics $1/\sqrt{c\left( y\right) }$ is the speed of sound. Let
the number $T>0.$ Consider the following Cauchy problem: 
\begin{align}
& c\left( y\right) u_{tt}=u_{yy},\quad \left( y,t\right) \in \mathbb{R}%
\times \left( 0,T\right) ,  \label{2.3} \\
& u\left( y,0\right) =0,\quad u_{t}\left( y,0\right) =\delta \left( y\right)
.  \label{2.4}
\end{align}%
The problem of finding the function $u\left( y,t\right) $ from conditions (%
\ref{2.3}), (\ref{2.4}) is our forward problem. Our interest is in the
following inverse problem:

\vspace{1em} \textbf{Coefficient Inverse Problem 1 (CIP1)}. \emph{Suppose
that the following two functions } \newline
$g_{0}\left( t\right) ,g_{1}\left( t\right) $\emph{\ are known:}%
\begin{equation}
u\left( 0,t\right) =g_{0}\left( t\right) ,\quad u_{y}\left( 0,t\right)
=g_{1}\left( t\right) ,\quad t\in \left( 0,T\right) .  \label{2.5}
\end{equation}
\emph{Determine the function }$c\left( y\right) $\emph{\ for }$y\in \left(
0,1\right) ,$ \emph{assuming that the number }$\overline{c}$\emph{\ is
known. }\vspace{1em}

We now introduce a change of variables in order to reduce the hyperbolic
equation (\ref{2.3}) to the wave-like equation with the unknown potential
function and the constant principal part of the 1D hyperbolic operator. This
is a well known change of variables, see, e.g. formulas (2.141)-(2.143) in
section 7 of chapter 2 of the book of Romanov \cite{romanov2018inverse}.
Thus, we introduce a new variable $x$ as 
\begin{equation}
x=\int \displaylimits_{0}^{y}\sqrt{c\left( s\right) }ds.  \label{2.6}
\end{equation}%
Physically, $x\left( y\right) $ is the travel time needed for the wave to
travel from the source $\left\{ 0\right\} $ to the point $\left\{ y\right\}
. $ Denote 
\begin{align}
& v\left( x,t\right) =u\left( y\left( x\right) ,t\right) c^{1/4}\left(
y\left( x\right) \right) ,  \label{2.81} \\
& S\left( x\right) =c^{-1/4}\left( y\left( x\right) \right) ,\quad r\left(
x\right) =\frac{S^{\prime \prime }\left( x\right) }{S\left( x\right) }-2%
\left[ \frac{S^{\prime }\left( x\right) }{S\left( x\right) }\right] ^{2}.
\label{2.8}
\end{align}%
Then, using (\ref{2.3}), (\ref{2.4}) and (\ref{2.6})-(\ref{2.8}), we obtain 
\begin{align}
& v_{tt}=v_{xx}+r\left( x\right) v,\quad \left( x,t\right) \in \mathbb{R}%
\times (0,\widetilde{T}),  \label{2.9} \\
& v\left( x,0\right) =0,\quad v_{t}\left( x,0\right) =\delta \left( x\right)
,  \label{2.10}
\end{align}%
where the number $\widetilde{T}=\widetilde{T}\left( T\right) $ depends on $T$%
. Using again (\ref{2.1}), (\ref{2.2}) and (\ref{2.5})-(\ref{2.8}), we
obtain 
\begin{align}
& r\left( x\right) =0\text{ for }x\in \left( -\infty ,0\right] \cup \left[
b,\infty \right) ,\quad r\left( x\right) \in C^{1}\left( \mathbb{R}\right) ,
\label{2.11} \\
& v\left( 0,t\right) =g_{0}\left( t\right) ,\quad v_{x}\left( 0,t\right)
=g_{1}\left( t\right) ,  \label{2.12} \\
& b=\int \displaylimits_{0}^{1}\sqrt{c\left( s\right) }ds.  \nonumber
\end{align}%
Existence and uniqueness of the solution of the forward problem (\ref{2.9})-(%
\ref{2.10}) is well known. More precisely, it was proven in, e.g. section 3
of chapter 2 of \cite{romanov2018inverse}, that 
\begin{align}
& v\left( x,t\right) =\frac{1}{2}+\frac{1}{2}\int \displaylimits_{\left(
x-t\right) /2}^{\left( x+t\right) /2}r\left( \xi \right) \int \displaylimits%
_{\left\vert \xi \right\vert }^{t-\left\vert x-\xi \right\vert }v\left( \xi
,\tau \right) d\tau ,t\geq \left\vert x\right\vert ,  \label{2.13} \\
& v\left( x,t\right) =0,\text{ }t<\left\vert x\right\vert ,  \nonumber \\
& \lim_{t\rightarrow \left\vert x\right\vert ^{+}}v\left( x,t\right) =\frac{1%
}{2}.  \label{2.14}
\end{align}%
Equation (\ref{2.13}) is a Volterra-like integral equation of the second
kind. Hence, it can be solved via sequential iterations. The corresponding
series converges absolutely in appropriate bounded subdomains of $\left\{
t\geq \left\vert x\right\vert \right\} $ together with its derivatives up to
the third order \cite{romanov2018inverse}. The existence of third
derivatives of iterates is guaranteed by $r\left( x\right) \in C^{1}\left( 
\mathbb{R}\right) $ and (\ref{2.13}). Thus, 
\begin{equation}
v\in C^{3}\left( t\geq \left\vert x\right\vert \right) .  \label{2.15}
\end{equation}%
If we would find the function $r\left( x\right) ,$ then the function $%
c\left( y\right) $ can be reconstructed via backward calculations, see
section 7 for details. Therefore, we arrive at CIP2.

\vspace{1em} \textbf{Coefficient Inverse Problem 2 (CIP2)}. \emph{Suppose
that the number }$\overline{c}$ in \emph{(\ref{2.1}) and the functions }$%
g_{0}\left( t\right) ,g_{1}\left( t\right) $\emph{\ in (\ref{2.12}) are
known. Determine the coefficient }$r\left( x\right) \in C^{1}\left( \mathbb{R%
}\right) $ \emph{in (\ref{2.9}) for }$x\in \left( 0,b\right) .$ \emph{\ }

\vspace{0.5em}

Note that the number $b$ is unknown and should be determined when solving
CIP2 and also (\ref{2.1}) implies that $b\geq 1.$ However, we can estimate $b
$ using (\ref{2.1}), (\ref{2.2}) as $b\leq \sqrt{\overline{c}}.$ Hence, we
fix a number $a$ such that 
\begin{equation}
a \geq \sqrt{\overline{c}} \geq b \geq 1.  \label{2.16}
\end{equation}%
We use this number $a$ everywhere below. It was established in \cite{romanov2018inverse} that $\widetilde{T}\geq 2b$ guarantees uniqueness of
CIP2. Hence, we take below 
\begin{equation}
\widetilde{T}=2a.  \label{2.17}
\end{equation}

Traditionally, classical absorbing boundary conditions of Engquist and Majda 
\cite{engquist1977absorbing} are artificially imposed when working with the
propagation of waves. In our case, however, Lemma 1 guarantees that the
absorbing boundary condition is satisfied indeed at any two point $x_{1}\geq
b$ and $x_{2}\leq 0$. We see in our computations that this condition is
important, since it provides a better stability property.

\textbf{Lemma 1} (absorbing boundary conditions). \emph{Let }$x_{1}\geq b$%
\emph{\ and }$x_{2}\leq 0$ \emph{be two arbitrary number.} \emph{Then the
solutions }$u\left( x,t\right) $\ \emph{of problem (\ref{2.3}), (\ref{2.4})
and }$v\left( x,t\right) $\emph{\ of problem (\ref{2.9}), (\ref{2.10})
satisfy the absorbing boundary condition at }$x=x_{1}$ and $x=x_{2}$\emph{\
i.e.}
\begin{align}
v_{x}(x_{1},t)+v_{t}(x_{1},t) &=0\text{ for }t\in (0,\widetilde{T}),
\label{2.18} \\
v_{x}\left( x_{2},t\right) -v_{t}\left( x_{2},t\right)  &=0\text{ for }t\in
(0,\widetilde{T}),  \label{2.19} \\
u_{x}\left( x_{1},t\right) +u_{t}\left( x_{1},t\right)  &=0\text{ for }t\in
(0,\widetilde{T}),  \label{2.20} \\
u_{x}\left( x_{2},t\right) -u_{t}\left( x_{2},t\right)  &=0\text{ for }t\in
(0,\widetilde{T}).  \label{2.21}
\end{align}

We omit the proof of Lemma 1 since it mainly follows from Lemma 2.2 of \cite{smirnov2020convexification}. Indeed, in the cases of (\ref{2.18}) and (\ref{2.19}) the proof of Lemma 1 is completely similar with the proof of Lemma 2.2 of \cite{smirnov2020convexification}. In the cases of (\ref{2.20}) and (\ref{2.21}), in the proof of Lemma 2.2 of \cite{smirnov2020convexification},
one should take into account (\ref{2.81}), (\ref{2.8}) and (\ref{2.11}), which can be done easily. 

\section{A Quasilinear PDE with a Non Local Term}

\label{sec:3}

Consider the rectangle $\Omega \subset \mathbb{R}^{2},$%
\begin{equation}
\Omega =\left( 0,a\right) \times (0,\widetilde{T}),\quad \widetilde{T}=2a.
\label{3.0}
\end{equation}%
Here, $\widetilde{T}$ is as in (\ref{2.17}). Consider the function $w$
defined as: 
\begin{equation}
w\left( x,t\right) =v\left( x,x+t\right) .  \label{3.1}
\end{equation}%
Then, using (\ref{2.9}) and (\ref{2.14}), we obtain%
\begin{equation}
w_{xx}-2w_{xt}=-r\left( x\right) w,\quad \left( x,t\right) \in \Omega
\label{3.2}
\end{equation}%
and $w\left( x,0\right) =1/2.$ Hence, 
\begin{equation}
r\left( x\right) =4w_{xt}\left( x,0\right) ,\text{ }x\in \left( 0,a\right) .
\label{3.3}
\end{equation}
Differentiate both sides of equation (\ref{3.2}) with respect to $t$, see (%
\ref{2.15}) and denote 
\begin{equation}
q\left( x,t\right) =w_{t}\left( x,t\right) ,\text{ }\left( x,t\right) \in
\Omega .  \label{3.30}
\end{equation}%
Then (\ref{2.5}), (\ref{2.12}), (\ref{2.14}) and (\ref{3.1})-(\ref{3.30})
lead to: 
\begin{align}
& q_{xx}-2q_{xt}+4q_{x}\left( x,0\right) q=0,\quad \left( x,t\right) \in
\Omega ,  \label{3.4} \\
& q\left( 0,t\right) =s_{0}\left( t\right) ,\quad q_{x}\left( 0,t\right)
=s_{1}\left( t\right) ,\quad q_{x}(a,t)=0,\quad t\in (0,\widetilde{T}),
\label{3.5} \\
& s_{0}\left( t\right) =g_{0}^{\prime }\left( t\right) ,s_{1}\left( t\right)
=g_{0}^{\prime \prime }\left( t\right) +g_{1}^{\prime }\left( t\right) .
\label{3.50}
\end{align}%
The condition $q_{x}\left( a,t\right) =0$ follows from (\ref{2.16}), Lemma
1 and (\ref{3.1}). The problem (\ref{3.4})-(\ref{3.5}) is a boundary value
problem (BVP) for a 1D quasilinear hyperbolic equation (\ref{3.4}) with
overdetermined boundary conditions (\ref{3.5}), an absent initial condition
at $\left\{ t=0\right\} $ and a non local term $q\left( x,0\right) .$

Suppose that we have found the function $q\left( x,t\right) $ satisfying (%
\ref{3.4})-(\ref{3.5}). Then, using (\ref{3.3}) and (\ref{3.30}), we obtain 
\begin{equation}
r\left( x\right) =4q_{x}\left( x,0\right) ,\text{ }x\in \left( 0,a\right) .
\label{3.6}
\end{equation}%
Therefore, we focus below on the numerical solution of the BVP (\ref{3.4})-(%
\ref{3.5}) via the convexification method.

\section{Globally Strictly Convex Tikhonov-like Cost Functional\ }

\label{sec:4}

\subsection{The main idea of the convexification}

\label{sec:4.1}

The first step of the convexification basically consists in obtaining an over-determined boundary value problem (BVP) for either a quasilinear integral differential equation with Volterra integrals in it \cite{KlibanovKamburg:mmas2016,klibanov2019backscatter,klibanovzhang2020convexification,smirnov2020convexification} or for a coupled system of elliptic PDEs \cite{khoa2019convexification,klibanov2019convexification,klibanov2018electrical,klibanov2018new}. The common point of both is that neither of these PDEs contains the
unknown coefficient. The latter reminds one the first step of the method of \cite{BukhgeimKlibanov:smd1981}.

Next, to solve this BVP, a weighted Tikhonov-like functional $S_{\lambda }$
is constructed, where $\lambda \geq 1$ is the parameter. The weight is the
Carleman Weight Function (CWF). This is the function, which is involved as
the weight in the Carleman estimate for the corresponding PDE operator. This
functional is evaluated on a convex bounded set $Q\left( 2d\right) \subset
H^{k}$ of the diameter $2d,$ where $d>0$ is an arbitrary but fixed number
and $H^{k},k\geq 1$ is an appropriate Sobolev space. The central theorem
then is that, for a proper choice $\lambda \left( d\right) $ of the
parameter of the functional $S_{\lambda },$ this functional is strictly
convex on $Q\left( 2d\right) $ for all $\lambda \geq \lambda \left( d\right)
.$ This eliminates the above mentioned phenomenon of multiple local minima.
Next, results of \cite{KlibanovNik:ra2017} enable one to prove existence and
uniqueness of the minimizer of $S_{\lambda }$ on $\overline{Q\left(
2d\right) }$ as well as the convergence of the gradient projection method to
the exact solution if starting from an arbitrary point of the set $Q\left(
2d\right) .$ Since the number $d$ is an arbitrary one, then this is the
global convergence, see section 1. We also note that even though the theory
requires the parameter $\lambda $ to be sufficiently large, our extensive
computational experience with the convexification, including the current
paper (see section 7.1), tells us that $\lambda \in \left[ 1,3\right] $
provides quite accurate reconstructions \cite{khoa2019convexification,klibanov2019backscatter,klibanovzhang2020convexification,klibanov2019convexification,klibanov2018electrical,klibanov2018new,smirnov2020convexification}.

\subsection{The Tikhonov-like Cost Functional for BVP (\protect\ref{3.4})-(%
\protect\ref{3.5})}

\label{sec:4.2}

\begin{figure}[tbh]
\begin{center}
\subfloat[The rectangle $R\left( 0,t\right)$.]{\includegraphics[width
=.45\textwidth]{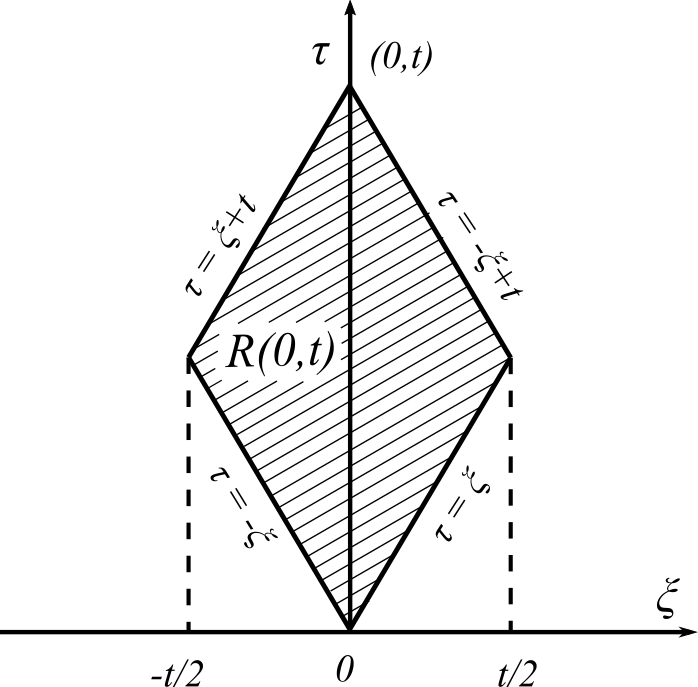}} \quad 
\subfloat[The triangle
$R_{a}$.]{\includegraphics[width =.38\textwidth]{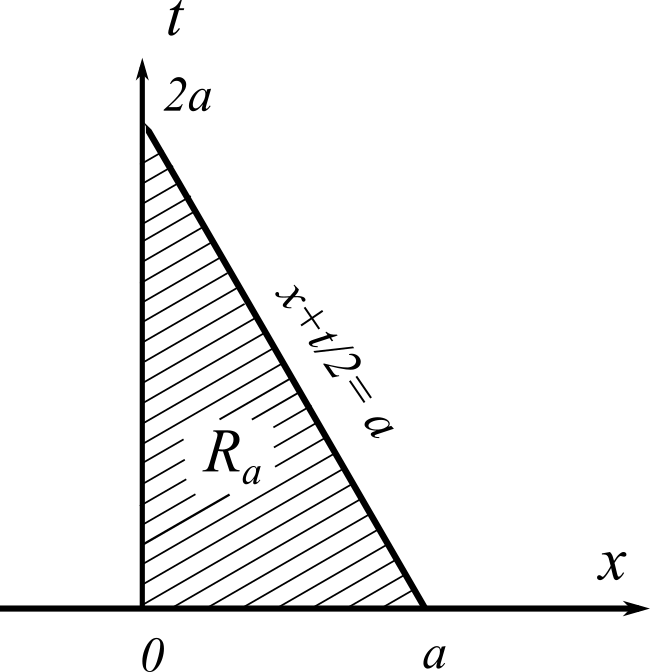}}
\end{center}
\caption{\emph{The rectangle $R\left( 0,t\right) $\emph{\ (see (\ref{4.1}))} and the triangle $R_{a}$\emph{\ (see (\protect\ref{4.2})).}}}
\label{fig1}
\end{figure}

Since $r\left( x\right) =0$ for $x<0,$ the integration in (\ref{2.13}) is
carried out over the rectangle 
\begin{equation}
R\left( x,t\right) =\left\{ \left( \xi ,\tau \right) :0<\xi <\tau
<t-\left\vert x-\xi \right\vert \right\} .  \label{4.1}
\end{equation}%
The change of variables (\ref{3.1}) transforms $R\left( 0,t\right) $ to the
following triangle \newline
$\left\{ \left( \xi ,\tau \right) :\xi ,\tau >0,\xi +\tau /2<t/2\right\} .$
Hence, it follows from (\ref{3.6}) that if we would find the function $%
q\left( x,t\right) $ in the triangle $R_{a}$ (see Figure \ref{fig1}), then
we would find the function $r\left( x\right) $ for $x\in \left( 0,b\right) $%
. Here, 
\begin{equation}
R_{a}=\left\{ \left( x,t\right) :x,t>0,x+t/2<a\right\} .  \label{4.2}
\end{equation}%

The CWF in this paper is the same as the one in \cite%
{smirnov2020convexification},%
\begin{equation}
\psi _{\lambda }\left( x,t\right) =e^{-2\lambda \left( x+\alpha t\right)
},\quad \alpha \in \left( 0,1/2\right) ,\quad \lambda \geq 1,\quad
\label{4.3}
\end{equation}%
where $\alpha $ and $\lambda $ are parameters independent on $x,t$.
Following (\ref{3.4}), consider the quasilinear operator $M$, 
\begin{equation}
M\left( q\right) =q_{xx}-2q_{xt}+4q_{x}\left( x,0\right) q,\text{ }\left(
x,t\right) \in \Omega .  \label{4.5}
\end{equation}%
Let $d>0$ be an arbitrary number. Consider the convex set $P\left(
d,s_{0},s_{1}\right) $ of the diameter $2d,$%
\begin{equation}
P\left( d,s_{0},s_{1}\right) =\left\{ 
\begin{array}{l}
q\in H^{4}\left( \Omega \right) ,\quad q\left( 0,t\right) =s_{0}\left(
t\right) ,\quad q_{x}\left( 0,t\right) =s_{1}\left( t\right) , \\ 
q_{x}\left( a,t\right) =0,\quad \left\Vert q\right\Vert _{H^{3}\left( \Omega
\right) }<d,%
\end{array}%
\right\}  \label{4.6}
\end{equation}%
where functions $s_{0}\left( t\right) ,s_{1}\left( t\right) $ are defined in
(\ref{3.50}). Our Tikhonov-like cost functional functional is%
\begin{equation}
K_{\lambda ,\gamma }\left( q\right) =\int \displaylimits_{\Omega }\left[
M\left( q\right) \right] ^{2}\psi _{\lambda }dxdt+\gamma \left\Vert
q\right\Vert _{H^{4}\left( \Omega \right) }^{2},  \label{4.7}
\end{equation}%
where $\gamma \in \left( 0,1\right) $ is the regularization parameter. The
reason why we use the space $H^{4}\left( \Omega \right) $ in the penalty
term of (\ref{4.7}) is that in the proof of Theorem 2 (below), we require \newline $P\left( d,s_{0},s_{1}\right) \subset C^{2}\left( \overline{\Omega }\right) $
as well as%
\begin{equation}
\left\Vert q\right\Vert _{C^{2}\left( \overline{\Omega }\right) }\leq
C_{1}d,\quad \forall q\in P\left( d,s_{0},s_{1}\right) .  \label{4.8}
\end{equation}%
Embedding theorem guarantees these. At the same time, we have established
numerically that we can use the discrete analog of $H^{2}\left( \Omega
\right) $ in the penalty term in our computations, see section 7. Here and
below $C_{1}=C_{1}\left( \Omega \right) >0$ denotes different positive
numbers depending only on the rectangle $\Omega .$ We concentrate below on
the solution of the following minimization problem:

\vspace{1em} \textbf{Minimization Problem}. \emph{Find a minimizer of the
functional }$K_{\lambda ,\gamma }\left( q\right) $\emph{\ on the set } 
\newline
$P\left( d,s_{0},s_{1}\right) .$ \vspace{1em}

We introduce two functional spaces $H_{0}^{2}\left( \Omega \right) \subset
H^{2}\left( \Omega \right) $ and $H_{0}^{3}\left( \Omega \right) \subset
H^{3}\left( \Omega \right) ,$ 
\begin{align}
& H_{0}^{2}\left( \Omega \right) =\left\{ u\in H^{2}\left( \Omega \right)
:u\left( 0,t\right) =u_{x}\left( 0,t\right) =0,\text{ }t\in (0,\widetilde{T}%
)\right\} ,  \nonumber \\
& H_{0}^{4}\left( \Omega \right) =\left\{ u\in H^{4}\left( \Omega \right)
:u\left( 0,t\right) =u_{x}\left( 0,t\right) =u_{x}\left( a,t\right) =0,\text{
}t\in (0,\widetilde{T})\right\} ,  \nonumber
\end{align}%
where the condition $u_{x}\left( a,t\right) =0$ is introduced due to the
last condition (\ref{3.5}).

\section{Convergence Analysis}

\label{sec:5}

In this section we formulate some theorems, which provide the convergence
analysis of our method of solving the above Minimization Problem. \vspace{1em%
}

\textbf{Theorem 1} (Carleman estimate). \emph{Let }$\psi _{\lambda }\left(
x,t\right) $\emph{\ be the function defined in (\ref{4.3}). Then there
exists a number }$C=C\left( \alpha ,\Omega \right) >0$\emph{\ and another
number } $\lambda _{0}=\lambda _{0}\left( \alpha ,\Omega \right) \geq 1$%
\emph{\ depending only on listed parameters},\emph{\ such that for all
functions }$w\in H_{0}^{2}\left( \Omega \right) $\emph{\ and for all }$%
\lambda \geq \lambda _{0}$\emph{\ the following Carleman estimate is valid:}%
\begin{equation*}
\begin{split}
& \int \displaylimits_{\Omega }\left( w_{xx}-2w_{xt}\right) ^{2}\psi
_{\lambda }dxdt\geq C\int \displaylimits_{\Omega }\left[ \lambda \left(
w_{x}^{2}+w_{t}^{2}\right) +\lambda ^{3}w^{2}\right] \psi _{\lambda }dxdt \\
& +C\int \displaylimits_{0}^{a}\left[ \lambda w_{x}^{2}+\lambda ^{3}w^{2}%
\right] \left( x,0\right) e^{-2\lambda x}dx-Ce^{-2\lambda \alpha \widetilde{T%
}}\int \displaylimits_{0}^{a}\left[ \lambda w_{x}^{2}+\lambda ^{3}w^{2}%
\right] \left( x,\widetilde{T}\right) dx.
\end{split}%
\end{equation*}

Let $R_{a}\subset \Omega $ be the triangle defined in (\ref{4.2}), where $%
\Omega $ is the rectangle defined in (\ref{3.0}).\emph{\ }Choose an
arbitrary number $\mu \in \left( 0,2\alpha a\right) .$ Define the triangle $%
R_{a,\alpha ,\mu }$ as%
\begin{equation}
R_{a,\alpha ,\mu }=\left\{ \left( x,t\right) :x+\alpha t<2\alpha a-\mu ,%
\hspace{0.3em}x,t>0\right\} \subset R_{a}.  \label{5.2}
\end{equation}%
Let $H$ be a Hilbert space, $\Phi \subset H$ be a convex set in $H$ and $%
J:\Phi \rightarrow H$ be a functional, which has the Frech\'{e}t derivative $%
J^{\prime }$ at any point $\mathbf{x}\in H.$ We remind that the following
estimate guarantees the strict convexity of $J$ on the set $\Phi :$%
\begin{equation}
J\left( \mathbf{y}\right) -J\left( \mathbf{x}\right) -J^{\prime }\left( 
\mathbf{x}\right) \left( \mathbf{y}-\mathbf{x}\right) \geq \beta \left\Vert 
\mathbf{x}-\mathbf{y}\right\Vert _{H}^{2},  \label{5.20}
\end{equation}%
for all $\mathbf{x},\mathbf{y}\in \Phi $ \cite{polyak1987}. Here, the
constant $\beta >0$ is independent on $\mathbf{x,y}$. \vspace{1em}

\textbf{Theorem 2} (the central theorem: global strict convexity).\emph{\
For any pair }$\lambda ,\gamma >0$\emph{\ and for any function }$q\in 
\overline{P\left( d,s_{0},s_{1}\right) }$\emph{\ the functional }$K_{\lambda
,\gamma }\left( q\right) $\emph{\ has the Frech\'{e}t derivative }$%
K_{\lambda ,\gamma }^{\prime }\left( q\right) \in H_{0}^{4}\left( \Omega
\right) .$\emph{\ Let } $\lambda _{0}=\lambda _{0}\left( \alpha \right) \geq
1$\emph{\ be the number of Theorem 1 and let the number }$\mu \in \left(
0,2\alpha a\right) .$\emph{\ There exist a sufficiently large number }$%
\lambda _{1}=\lambda _{1}\left( \alpha ,d,\Omega \right) \geq $\emph{\ }$%
\lambda _{0}$\emph{\ and a number }$B=B\left( \alpha ,d,\Omega \right) >0$%
\emph{, both depending only on listed parameters, such that for all }$%
\lambda \geq \lambda _{1}$\emph{\ and for all } \newline $\gamma \in \lbrack
2e^{-\lambda \alpha \widetilde{T}},1),$\emph{\ the functional }$K_{\lambda
,\gamma }\left( q\right) $\emph{\ is strictly convex on the set }$\overline{%
P\left( d,s_{0},s_{1}\right) }$\emph{, i.e. (see (\ref{5.20}))\ }%
\begin{equation*}
\begin{split}
& K_{\lambda ,\gamma }\left( q_{2}\right) -K_{\lambda ,\gamma }\left(
q_{1}\right) -K_{\lambda ,\gamma }^{\prime }\left( q_{1}\right) \left(
q_{2}-q_{1}\right) \geq Be^{-2\lambda \left( 2\alpha a-\mu \right)
}\left\Vert q_{2}-q_{1}\right\Vert _{H^{1}\left( R_{a,\alpha ,\mu }\right)
}^{2} \\
& +Be^{-2\lambda \left( 2\alpha a-\mu \right) }\left\Vert q_{2}\left(
x,0\right) -q_{1}\left( x,0\right) \right\Vert _{H^{1}\left( 0,2\alpha a-\mu
\right) }^{2} \\
& +\frac{\gamma }{2}\left\Vert q_{2}-q_{1}\right\Vert _{H^{4}\left( \Omega
\right) }^{2},\text{ }\forall q_{1},q_{2}\in \overline{P\left(
d,s_{0},s_{1}\right) },\quad \forall \lambda \geq \lambda _{1}.
\end{split}
\end{equation*}

Below $B=B\left( \alpha ,d,\Omega \right) >0$ and $C=C\left( \alpha ,\Omega
\right) >0$ denote different numbers depending only on listed parameters.

\vspace{1em} \textbf{\ Theorem 3}. \emph{Let parameters }$\lambda
_{1},\lambda ,\gamma $\emph{\ be the same as in Theorem 2. Then there exists
unique minimizer }$q_{\min ,\lambda ,\gamma }\in \overline{P\left(
d,s_{0},s_{1}\right) }$\emph{\ of the functional }$K_{\lambda ,\gamma
}\left( q\right) $\emph{\ on the set }$\overline{P\left(
d,s_{0},s_{1}\right) }.$\emph{\ Furthermore, the following inequality holds}%
\begin{equation*}
K_{\lambda ,\gamma }^{\prime }\left( q_{\min ,\lambda ,\gamma }\right)
\left( q-q_{\min ,\lambda ,\gamma }\right) \geq 0,\quad \forall q\in 
\overline{P\left( d,s_{0},s_{1}\right) }.
\end{equation*}

We now want to estimate the accuracy of the reconstruction in the presence
of noise in the data. Following one of the main concepts of the
regularization theory, we assume is the existence of the exact solution $%
r^{\ast }\left( x\right) \in C^{1}\left( \mathbb{R}\right) $ of the CIP2
with the ideal noiseless data \cite{BeilinaKlibanovBook,TikhonovGoncharsky}. Furthermore, we
assume that conditions (\ref{2.11}) hold for the function $r^{\ast }\left(
x\right) $. Let $q^{\ast }\left( x,t\right) $ be the corresponding function $%
q\left( x,t\right) $. We assume that $q^{\ast }\in P\left( d,s_{0}^{\ast
},s_{1}^{\ast }\right) ,$ where functions $s_{0}^{\ast },s_{1}^{\ast }$ are
the noiseless data $s_{0},s_{1}$ in (\ref{3.5}). It follows from (\ref{3.4})
and (\ref{3.6}) that 
\begin{equation}
M\left( q^{\ast }\right) =0,\hspace{0.3em}\left( x,t\right) \in \Omega ,%
\text{ and }r^{\ast }\left( x\right) =4q_{x}^{\ast }\left( x,0\right) ,\text{
}x\in \left( 0,a\right) .  \label{5.301}
\end{equation}%
Let $\sigma $ be the level of noise in the data. We assume that the number $%
\sigma ,$ 
\begin{equation}
\sigma \in \left( 0,\min \left( d,1\right) \right)   \label{5.30}
\end{equation}%
is the noise level in the data. Suppose that there exist functions $Q\in
P\left( d,s_{0},s_{1}\right) $ and \newline
$Q^{\ast }\in P\left( d,s_{0}^{\ast },s_{1}^{\ast }\right) $ such that 
\begin{equation}
\left\Vert Q-Q^{\ast }\right\Vert _{H^{4}(R_{a})}<\sigma .  \label{5.3}
\end{equation}%
Introduce functions $Y^{\ast }$ and $Y$ as 
\begin{align*}
& Y^{\ast }=q^{\ast }-Q^{\ast },  \label{5.4} \\
& Y=q-Q. 
\end{align*}%
Let $D>0$ be an arbitrary number. Define 
\begin{equation}
P_{0}\left( D\right) =\left\{ w\in H_{0}^{4}\left( \Omega \right)
:\left\Vert w\right\Vert _{H^{4}\left( \Omega \right) }<D\right\} .
\label{5.6}
\end{equation}%
Using (\ref{5.30})-(\ref{5.6}) and the triangle inequality, we obtain%
\begin{align}
& Y^{\ast }\in P_{0}\left( 2d\right) ,  \nonumber \\
& Y\in P_{0}\left( 2d\right) ,\text{ }\forall q\in P\left(
d,s_{0},s_{1}\right) ,  \label{5.80} \\
& Y+Q\in P\left( 3d,s_{0},s_{1}\right) ,\text{ }\forall Y\in P_{0}\left(
2d\right) .  \label{5.90}
\end{align}%
We now consider a modification $\widetilde{K}_{\lambda ,\gamma }:P_{0}\left(
2d\right) \rightarrow \mathbb{R}$ of the functional $K_{\lambda ,\gamma },$ 
\begin{equation}
\widetilde{K}_{\lambda ,\gamma }\left( Y\right) =K_{\lambda ,\gamma }\left(
Y+Q\right) ,\quad \forall Y\in P_{0}\left( 2d\right) .  \label{5.10}
\end{equation}

\textbf{Theorem 4}. \emph{The functional }$\widetilde{K}_{\lambda ,\gamma }$%
\emph{\ has the Frech\'{e}t derivative }$\widetilde{K}_{\lambda ,\gamma
}\left( Y\right) \in H_{0}^{3}(\Omega )$\emph{\ for every point }$Y\in 
\overline{P_{0}\left( 2d\right) }$\emph{\ and for all }$\lambda ,\gamma >0.$%
\emph{\ Denote }$\lambda _{2}=\lambda _{1}\left( \alpha ,3d,\Omega \right)
\geq \lambda _{1},$\emph{\ where }$\lambda _{1}=\lambda _{1}\left( \alpha
,d,\Omega \right) \geq 1$\emph{\ is the number of Theorem 2.\ For every }$%
\lambda \geq \lambda _{2}$\emph{\ and for every }$\gamma \in \lbrack
2e^{-\lambda \alpha \widetilde{T}},1)$\emph{\ the functional }$\widetilde{K}%
_{\lambda ,\gamma }\left( Y\right) $\emph{\ is strictly convex on the ball }$%
\overline{P_{0}\left( 2d\right) }\subset H_{0}^{4}(\Omega ),$ i.e. \emph{%
(see (\ref{5.20})),} 
\begin{equation}
\begin{split}
& \widetilde{K}_{\lambda ,\gamma }\left( Y_{2}\right) -\widetilde{K}%
_{\lambda ,\gamma }\left( Y_{1}\right) -\widetilde{K}_{\lambda ,\gamma
}^{\prime }\left( Y_{1}\right) \left( Y_{2}-Y_{1}\right) \geq Be^{-2\lambda
\left( 2\alpha b-\mu \right) }\left\Vert Y_{2}-Y_{1}\right\Vert
_{H^{1}\left( R_{a,\alpha ,\mu }\right) }^{2} \\
& +Be^{-2\lambda \left( 2\alpha a-\mu \right) }\left\Vert Y_{2}\left(
x,0\right) -Y_{1}\left( x,0\right) \right\Vert _{H^{1}\left( 0,2\alpha a-\mu
\right) }^{2} \\
& +\frac{\gamma }{2}\left\Vert Y_{2}-Y_{1}\right\Vert _{H^{4}\left( \Omega
\right) }^{2},\text{ }\forall Y_{1},Y_{2}\in \overline{P_{0}\left( 2d\right) 
},\quad \forall \lambda \geq \lambda _{2}.
\end{split}
\label{500}
\end{equation}%
\emph{Furthermore, there exists unique minimizer }$Y_{\min ,\lambda ,\gamma
} $\emph{\ }$\in \overline{P_{0}\left( 2d\right) }$\emph{\ of the functional 
}$\widetilde{K}_{\lambda ,\gamma }\left( Y\right) $\emph{\ on the set }$%
\overline{P_{0}\left( 2d\right) },$ \emph{and the following inequality holds}%
\begin{equation}
\widetilde{K}_{\lambda ,\gamma }^{\prime }\left( Y_{\min ,\lambda ,\gamma
}\right) \left( Y-Y_{\min ,\lambda ,\gamma }\right) \geq 0,\quad \forall
Y\in \overline{P_{0}\left( 2d\right) }.  \label{5.11}
\end{equation}

\textbf{Theorem 5 }(accuracy estimate).\emph{\ Suppose that }$\widetilde{T}%
\geq 4a.$\emph{\ Denote }%
\begin{equation}
\nu =\frac{\alpha \left( \widetilde{T}-4a\right) +\mu }{2\left( 2\alpha
a-\mu \right) },\quad \kappa =\frac{1}{2}\min \left( \nu ,1\right) .
\label{5.12}
\end{equation}%
\emph{Let the number }$\sigma _{0}\in \left( 0,1\right) $\emph{\ be so small
that }$\ln \sigma _{0}^{-1/\left( 2\left( 2\alpha a-\mu \right) \right)
}\geq \lambda _{2},$\emph{\ where }$\lambda _{2}$\emph{\ is the number of} 
\emph{Theorem 4.} \emph{Let }$\sigma \in \left( 0,\sigma _{0}\right) .$\emph{%
\ Let the \ the numbers }$\lambda =\lambda \left( \sigma \right) $ \emph{and}
$\gamma =\gamma \left( \sigma \right) $ \emph{be such that} \emph{\ }%
\begin{align}
& \lambda =\lambda \left( \sigma \right) =\ln \sigma ^{-1/\left( 2\left(
2\alpha a-\mu \right) \right) }>\lambda _{2},  \label{5.13} \\
& \gamma =\gamma \left( \sigma \right) =2e^{-\lambda \alpha \widetilde{T}%
}=2\sigma ^{\left( \alpha \widetilde{T}\right) /\left( 2\left( 2\alpha a-\mu
\right) \right) }.  \label{5.14}
\end{align}%
\emph{Let} $Y_{\min ,\lambda ,\gamma }\in \overline{P_{0}\left( 2d\right) }$ 
\emph{be the minimizer of the functional }$\widetilde{K}_{\lambda ,\gamma
}\left( Y\right) $\emph{\ on the set }$\overline{P_{0}\left( 2d\right) },$%
\emph{\ the existence and uniqueness of which is guaranteed by Theorem 4.} 
\emph{Denote}%
\begin{align}
& q_{\min ,\lambda ,\gamma }=\left( Y_{\min ,\lambda ,\gamma }+Q\right) \in 
\overline{P\left( 3d,s_{0},s_{1}\right) },  \label{5.15} \\
& r_{\min ,\lambda ,\gamma }\left( x\right) =4\partial _{x}\left[ q_{\min
,\lambda ,\gamma }\left( x,0\right) \right] .  \label{5.16}
\end{align}%
\emph{Then the following estimates are valid }%
\begin{align}
& \left\Vert q_{\min ,\lambda ,\gamma }-q^{\ast }\right\Vert _{H^{1}\left(
R_{a,\alpha ,\mu }\right) }\leq B\sigma ^{\kappa },  \label{5.17} \\
& \left\Vert r_{\min ,\lambda ,\gamma }-r^{\ast }\right\Vert _{L_{2}\left(
0,2\alpha a-\mu \right) }\leq B\sigma ^{\kappa }.  \label{5.19}
\end{align}

To minimize the functional $\widetilde{K}_{\lambda ,\gamma }\left( Y\right) $
on the set $\overline{P_{0}\left( 2d\right) }\subset H_{0}^{4}\left( \Omega
\right) $ via the gradient projection mehod, we first consider the
orthogonal projection operator $Z:H_{0}^{4}\left( \Omega \right) \rightarrow 
\overline{P_{0}\left( 2d\right) }$ \ of the space $H_{0}^{4}\left( \Omega
\right) $ on the closed ball $\overline{P_{0}\left( 2d\right) }.$ Let $%
Y_{0}\in P_{0}\left( 2d\right) $ be an arbitrary point and the number $%
\omega \in \left( 0,1\right) .$ The sequence of the gradient projection
method is \cite{KlibanovNik:ra2017}:%
\begin{equation}
Y_{n}=Z\left( Y_{n-1}-\omega \widetilde{K}_{\lambda ,\gamma }\left(
Y_{n-1}\right) \right) ,\quad n=1,2,\dots  \label{5.21}
\end{equation}

\textbf{Theorem 6}. \emph{Let }$\lambda _{1}=\lambda _{1}\left( \alpha
,d,\Omega \right) \geq 1$\emph{\ be the number of Theorem 2 and let }$%
\lambda _{2}$ be \emph{the number of Theorem 4. Let }$\ $\emph{the numbers }$%
\kappa ,\sigma _{0}$ \emph{be the same as one in Theorem 5. Let }$\sigma \in
\left( 0,\sigma _{0}\right) $ \emph{and let the numbers }$\lambda =\lambda
\left( \sigma \right) $\emph{\ and }$\gamma =\gamma \left( \sigma \right) $%
\emph{\ be the same as in (\ref{5.13}) and (\ref{5.14}) respectively. Let }$%
Y_{\min ,\lambda ,\gamma }\in \overline{P_{0}\left( 2d\right) }$\emph{\ be
the minimizer of the functional }$\widetilde{K}_{\lambda ,\gamma }\left(
Y\right) $ \emph{on the set }$\overline{P_{0}\left( 2d\right) }$,\emph{\ the
existence and uniqueness of which} \emph{is guaranteed by Theorem 4. Let the
function }$q_{\min ,\lambda ,\gamma }\in \overline{P\left(
3d,s_{0},s_{1}\right) }$\emph{\ be the one defined in (\ref{5.15}). Consider
functions }$q_{n}=Y_{n}+Q\in \overline{P\left( 3d,s_{0},s_{1}\right) },$%
\emph{\ }$n=0,1,\dots ,$ \emph{see (\ref{5.80}) and (\ref{5.90}).\ Also, let 
}$r_{n}\left( x\right) $\emph{\ and }$r_{\min ,\lambda ,\gamma }\left(
x\right) $\emph{\ be the coefficients }$r\left( x\right) ,$\emph{\ which are
found from the functions }$q_{n}$\emph{\ and }$q_{\min ,\lambda ,\gamma }$%
\emph{\ via (\ref{3.6}) and (\ref{5.16}) respectively. Then there exists a
number }$\omega _{0}=\omega _{0}\left( \alpha ,\mu ,d,\sigma \right) \in
\left( 0,1\right) $\emph{\ depending only on listed parameters such that for
any }$\omega \in \left( 0,\omega _{0}\right) $\emph{\ there exists a number }%
$\theta =\theta \left( \omega \right) \in \left( 0,1\right) $\emph{\ such
that the following convergence rates are valid:}%
\begin{align}
& \left\Vert q_{\min ,\lambda ,\gamma }-q_{n}\right\Vert _{H^{3}\left(
\Omega \right) }\leq \theta ^{n}\left\Vert q_{\min ,\lambda ,\gamma
}-q_{0}\right\Vert _{H^{4}\left( \Omega \right) },\text{ }n=1,2,\dots ,
\label{5.22} \\
& \left\Vert r_{\min ,\lambda ,\gamma }\left( x\right) -r_{n}\right\Vert
_{L_{2\left( 0,2\alpha a-\mu \right) }}\leq \theta ^{n}\left\Vert q_{\min
,\lambda ,\gamma }-q_{0}\right\Vert _{H^{4}\left( \Omega \right) },\text{ }%
n=1,2,\dots ,  \label{5.23} \\
& \left\Vert q^{\ast }-q_{n}\right\Vert _{H^{1}\left( R_{a,\alpha ,\mu
}\right) }\leq B\sigma ^{\kappa }+\theta ^{n}\left\Vert q_{\min ,\lambda
,\gamma }-q_{0}\right\Vert _{H^{4}\left( \Omega \right) },\text{ }%
n=1,2,\dots ,  \label{5.24} \\
& \left\Vert r^{\ast }-r_{n}\right\Vert _{L_{2\left( 0,2\alpha a-\mu \right)
}}\leq B\sigma ^{\kappa }+\theta ^{n}\left\Vert q_{\min ,\lambda ,\gamma
}-q_{0}\right\Vert _{H^{4}\left( \Omega \right) },\text{ }n=1,2,\dots
\label{5.25}
\end{align}

\begin{remark}

\emph{Estimates (\ref{5.22})-(\ref{5.25}) imply that the sequence }$%
\left\{ Y_{n}\right\} $\emph{\ in (\ref{5.21}) generates the sequence of
coefficients }$\left\{ r_{n}\right\} ,$\emph{\ which converges globally to
the function }$r_{\min ,\lambda ,\gamma }$\emph{. Also, the sequence }$%
\left\{ r_{n}\right\} $\emph{\ converges globally to the exact coefficient }$%
r^{\ast }$\emph{\ as long as the level of the noise }$\xi $\emph{\ in the
data tends to zero. The global convergence property is due to the fact that
the starting point }$Y_{0}$\emph{\ of iterations (\ref{5.21}) is an
arbitrary point of the ball }$P_{0}\left( 2d\right) \subset H_{0}^{3}\left(
\Omega \right) $\emph{\ and the radius }$d$\emph{\ of this ball is an
arbitrary number.}
\end{remark}

\begin{remark}
\emph{Theorem 1 was proven in \cite{smirnov2020convexification}. Theorem 3 is a straightforward consequence of Theorem 2 combined with Lemma 2.1 of \cite{KlibanovNik:ra2017}. Moreover, it is obvious that Theorem 4 follows immediately from Theorems 2, 3 and (\ref{5.10}). Thus, we do not prove Theorems 1, 3 and 4 in this paper.}
\end{remark}

\begin{remark}
\emph{The most important difference between proofs of Theorems 2,5,6 and some theorems of \cite{smirnov2020convexification} is that, unlike \cite{smirnov2020convexification}, we do not work
here with an integral differential equation.}
\end{remark}

\section{Proofs of Theorems 2,5,6}

\label{sec:6}

In this section, $\left( x,t\right) \in \Omega ,$ where the rectangle $%
\Omega $ is defined in (\ref{3.0}).

\subsection{Proof of Theorem 2}

\label{sec:6.2}

Consider two arbitrary functions $q_{1},q_{2}\in \overline{P\left(
d,s_{0},s_{1}\right) }$. Let $h=q_{2}-q_{1}.$ By (\ref{4.6}), (\ref{5.6})
and the triangle inequality, 
\begin{equation}
h\in \overline{P_{0}\left( 2d\right) }.  \label{6.13}
\end{equation}%
As it was noticed in section 4, embedding theorem implies $\overline{P\left(
d,p_{0},p_{1}\right) },\overline{P_{0}\left( 2d\right) }\subset C^{2}\left( 
\overline{\Omega }\right) .$ Also, (\ref{4.8}) and (\ref{6.13}) lead to: 
\begin{equation}
\left\Vert q\right\Vert _{C^{2}\left( \overline{\Omega }\right) }\leq
B,\quad \forall q\in \overline{P\left( d,s_{0},s_{1}\right) },\quad
\left\Vert h\right\Vert _{C^{2}\left( \overline{\Omega }\right) }\leq B.
\label{6.14}
\end{equation}

First, we evaluate the expression $\left[ M\left( q_{1}+h\right) \right]
^{2}-\left[ M\left( q_{1}\right) \right] ^{2},$ 
\begin{align*}
& M\left( q_{1}+h\right) =\left( q_{1xx}-2q_{1xt}+4q_{1x}\left( x,0\right)
q_{1}\right) +\left( h_{xx}-2h_{xt}+4h_{x}\left( x,0\right)
q_{1}+4q_{1x}\left( x,0\right) h\right)  \\
& +4h_{x}\left( x,0\right) h=M\left( q_{1}\right) +\left(
h_{xx}-2h_{xt}+4h_{x}\left( x,0\right) q_{1}+4q_{1}\left( x,0\right)
h\right) +4h_{x}\left( x,0\right) h.
\end{align*}%
where the operator $M$ is defined in (\ref{4.5}). Hence,%
\begin{equation}
\begin{split}
& \left[ M\left( q_{1}+h\right) \right] ^{2}-\left[ M\left( q_{1}\right) %
\right] ^{2}=2M\left( q_{1}\right) \left( h_{xx}-2h_{xt}+4h_{x}\left(
x,0\right) q_{1}+4q_{1}\left( x,0\right) h\right)  \\
& +8M\left( q_{1}\right) h_{x}\left( x,0\right) h+\left[ \left(
h_{xx}-2h_{xt}+4h_{x}\left( x,0\right) q_{1}+4q_{1x}\left( x,0\right)
h\right) +4h_{x}\left( x,0\right) h\right] ^{2}.
\end{split}
\label{6.160}
\end{equation}%
In the right hand side of (\ref{6.160}), the first term is linear with
respect to $h$ and other terms are nonlinear respect to $h$. Let $%
M_{lin}\left( q_{1}\right) \left( h\right) $ be the linear part of the right
hand side of (\ref{6.160}), i.e. 
\begin{equation*}
M_{lin}\left( q_{1}\right) \left( h\right) =2M\left( q_{1}\right) \left(
h_{xx}-2h_{xt}+4h_{x}\left( x,0\right) q_{1}+4q_{1x}\left( x,0\right)
h\right) .  \label{6.17}
\end{equation*}%
Thus, 
\begin{equation}
\begin{split}
& \left[ M\left( q_{1}+h\right) \right] ^{2}-\left[ M\left( q_{1}\right) %
\right] ^{2}-M_{lin}\left( q_{1}\right) \left( h\right)  \\
& =8M\left( q_{1}\right) h_{x}\left( x,0\right) h+\left[ \left(
h_{xx}-2h_{xt}+4h_{x}\left( x,0\right) w_{1}+4q_{1x}\left( x,0\right)
h\right) +4h_{x}\left( x,0\right) h\right] ^{2}.
\end{split}
\label{6.180}
\end{equation}%
By (\ref{4.5}) and (\ref{6.14}) $\left\vert M\left( q_{1}\right) \right\vert
\leq B.$ Hence, the Cauchy-Schwarz inequality and (\ref{6.14}) imply that
the right hand side of (\ref{6.180}) can be estimated from the below as:%
\begin{equation*}
\begin{split}
\big[(h_{xx}-2h_{xt}& +4h_{x}(x,0)q_{1}+4q_{1x}(x,0)h)+4h_{x}(x,0)h\big]^{2}
\\
& +8M\left( q_{1}\right) h_{x}\left( x,0\right) h\geq \frac{1}{2}\left(
h_{xx}-2h_{xt}\right) ^{2}-B\left( h_{x}^{2}\left( x,0\right) +h^{2}\right) .
\end{split}
\label{6.190}
\end{equation*}%
It follows from (\ref{4.7}) and (\ref{6.180}) that%
\begin{equation}
\begin{split}
& K_{\lambda ,\gamma }\left( q_{1}+h\right) -K_{\lambda ,\gamma }\left(
q_{1}\right) =\int \displaylimits_{\Omega }M_{lin}\left( q_{1}\right) \left(
h\right) \psi _{\lambda }dxdt+2\left[ q_{1},h\right]  \\
& +\int \displaylimits_{\Omega }\left[ \left( h_{xx}-2h_{xt}+4h_{x}\left(
x,0\right) q_{1}+4q_{1x}\left( x,0\right) h\right) +4h_{x}\left( x,0\right) h%
\right] ^{2}\psi _{\lambda }dxdt \\
& +8\int \displaylimits_{\Omega }\left( M\left( q_{1}\right) h_{x}\left(
x,0\right) h\right) \psi _{\lambda }dxdt+\gamma \left\Vert h\right\Vert
_{H^{4}\left( \Omega \right) }^{2},
\end{split}
\label{6.191}
\end{equation}%
where $\left[ \cdot ,\cdot \right] $ denotes the scalar product in $%
H^{4}\left( \Omega \right) .$ Define the functional \newline
$A\left( q_{1}\right) :H_{0}^{4}\left( \Omega \right) \rightarrow \mathbb{R}$
as%
\begin{equation}
A\left( q_{1}\right) \left( z\right) =\int \displaylimits_{\Omega
}M_{lin}\left( q_{1}\right) \left( z\right) \psi _{\lambda }dxdt+2\gamma %
\left[ q_{1},z\right] ,\quad \forall z\in H_{0}^{4}\left( \Omega \right) .
\label{6.21}
\end{equation}%
Then $A\left( q_{1}\right) $ is a bounded linear functional. Next, it
follows from (\ref{6.191}) that 
\begin{equation*}
\lim_{\left\Vert y\right\Vert _{H^{4}\left( \Omega \right) }\rightarrow 0}%
\frac{1}{\left\Vert y\right\Vert _{H^{4}\left( \Omega \right) }}\left[
K_{\lambda ,\gamma }\left( q_{1}+y\right) -K_{\lambda ,\gamma }\left(
q_{1}\right) -A\left( q_{1}\right) \left( y\right) \right] =0.
\end{equation*}%
Hence, $A\left( q_{1}\right) $ is the Frech\'{e}t derivative of the
functional $K_{\lambda ,\gamma }\left( q\right) $ at the point $q_{1}.$ By
the Riesz theorem, there exists unique element 
\begin{equation*}
\widetilde{A}\left( q_{1}\right) \in H_{0}^{4}\left( \Omega \right) \hspace{%
0.3em}:\hspace{0.3em}A\left( q_{1}\right) \left( z\right) =\left[ \widetilde{%
A}\left( q_{1}\right) ,z\right] ,\quad \forall z\in H_{0}^{4}\left( \Omega
\right) .
\end{equation*}%
Thus, we set $\widetilde{A}\left( q_{1}\right) =K_{\lambda ,\gamma }^{\prime
}\left( q_{1}\right) \in H_{0}^{4}\left( \Omega \right) .$ Hence, using (\ref%
{6.190})-(\ref{6.21}), we obtain%
\begin{equation}
\begin{split}
& K_{\lambda ,\gamma }\left( q_{1}+h\right) -K_{\lambda ,\gamma }\left(
q_{1}\right) -K_{\lambda ,\gamma }^{\prime }\left( q_{1}\right) \left(
h\right) \geq \frac{1}{2}\int \displaylimits_{\Omega }\left(
h_{xx}-2h_{xt}\right) ^{2}\psi _{\lambda }dxdt \\
& -B\int \displaylimits_{\Omega }\left( h_{x}^{2}\left( x,0\right)
+h^{2}\right) \psi _{\lambda }dxdt+\gamma \left\Vert h\right\Vert
_{H^{4}\left( \Omega \right) }^{2}.
\end{split}
\label{6.210}
\end{equation}%
We now apply Carleman estimate of Theorem 1 to the right hand side of (\ref%
{6.210}),

\begin{equation}
\begin{split}
& \frac{1}{2}\int \displaylimits_{\Omega }\left( h_{xx}-2h_{xt}\right)
^{2}\psi _{\lambda }dxdt-B\int \displaylimits_{\Omega }\left(
h_{x}^{2}\left( x,0\right) +h^{2}\right) \psi _{\lambda }dxdt+\gamma
\left\Vert h\right\Vert _{H^{4}\left( \Omega \right) }^{2} \\
& \geq C\int \displaylimits_{\Omega }\left[ \lambda \left(
h_{x}^{2}+h_{t}^{2}\right) +\lambda ^{3}h^{2}\right] \psi _{\lambda
}dxdt-B\int \displaylimits_{\Omega }h^{2}\psi _{\lambda }dxdt \\
& +C\int \displaylimits_{0}^{a}\left[ \lambda h_{x}^{2}+\lambda ^{3}h^{2}%
\right] \left( x,0\right) e^{-2\lambda x}dx-B\int \displaylimits_{\Omega
}h_{x}^{2}\left( x,0\right) \psi _{\lambda }dxdt+\gamma \left\Vert
h\right\Vert _{H^{3}\left( \Omega \right) }^{2} \\
& -Ce^{-2\lambda \alpha \widetilde{T}}\int \displaylimits_{0}^{a}\left[
\lambda h_{x}^{2}+\lambda ^{3}h^{2}\right] \left( x,\widetilde{T}\right)
dx+\gamma \left\Vert h\right\Vert _{H^{4}\left( \Omega \right) }^{2}.
\end{split}
\label{6.23}
\end{equation}%
Next, 
\begin{equation*}
\int \displaylimits_{\Omega }h_{x}^{2}\left( x,0\right) \psi _{\lambda
}dxdt=\int \displaylimits_{0}^{a}h_{x}^{2}\left( x,0\right) e^{-2\lambda
x}dx\left( \int \displaylimits_{0}^{\widetilde{T}}e^{-2\lambda \alpha
t}dt\right) \leq \frac{1}{2\lambda \alpha }\int \displaylimits%
_{0}^{a}h_{x}^{2}\left( x,0\right) e^{-2\lambda x}dx.
\end{equation*}%
Hence, (\ref{6.23}) implies that%
\begin{equation}
\begin{split}
& \frac{1}{2}\int \displaylimits_{\Omega }\left( h_{xx}-2h_{xt}\right)
^{2}\psi _{\lambda }dxdt-B\int \displaylimits_{\Omega }\left(
h_{x}^{2}\left( x,0\right) +h^{2}\right) \psi _{\lambda }dxdt \\
& +\gamma \left\Vert h\right\Vert _{H^{4}\left( \Omega \right) }^{2}\geq
C\int \displaylimits_{\Omega }\left[ \lambda \left(
h_{x}^{2}+h_{t}^{2}\right) +\lambda ^{3}h^{2}\right] \psi _{\lambda
}dxdt-B\int \displaylimits_{\Omega }h^{2}\psi _{\lambda }dxdt \\
& +C\int \displaylimits_{0}^{a}\left[ \lambda h_{x}^{2}+\lambda ^{3}h^{2}%
\right] \left( x,0\right) e^{-2\lambda x}dx-\frac{B}{2\lambda \alpha }\int %
\displaylimits_{0}^{a}h_{x}^{2}\left( x,0\right) e^{-2\lambda x}dx+\gamma
\left\Vert h\right\Vert _{H^{4}\left( \Omega \right) }^{2} \\
& -C\lambda ^{3}e^{-2\lambda \alpha \widetilde{T}}\Vert h(x,\widetilde{T}%
)\Vert _{H^{1}\left( 0,a\right) }^{2}.
\end{split}
\label{6.24}
\end{equation}%
By the trace theorem, 
\begin{equation*}
\Vert z(x,\widetilde{T})\Vert _{H^{1}\left( 0,a\right) }^{2}\leq
C_{1}\left\Vert z\right\Vert _{H^{4}\left( \Omega \right) }^{2},\quad
\forall z\in H^{4}\left( \Omega \right) .  \label{6.25}
\end{equation*}%
Choose $\lambda _{1}=\lambda _{1}\left( \alpha ,d,\Omega \right) \geq
\lambda _{0}\geq 1$ so large that 
\begin{equation*}
e^{-\lambda \alpha \widetilde{T}}\geq CC_{1}\lambda ^{3}e^{-2\lambda \alpha 
\widetilde{T}},\quad C\lambda ^{3}\geq \frac{B}{\lambda \alpha }+2B,\quad
\forall \lambda \geq \lambda _{1}.
\end{equation*}%
Also, choose $\gamma \in \lbrack 2e^{-\lambda \alpha \widetilde{T}},1).$
Then (\ref{6.24}) implies that 
\begin{equation}
\begin{split}
& \frac{1}{2}\int \displaylimits_{\Omega }\left( h_{xx}-2h_{xt}\right)
^{2}\psi _{\lambda }dxdt-B\int \displaylimits_{\Omega }\left(
h_{x}^{2}\left( x,0\right) +h^{2}\right) \psi _{\lambda }dxdt \\
& +\gamma \left\Vert h\right\Vert _{H^{4}\left( \Omega \right) }^{2}\geq
B\int \displaylimits_{\Omega }\left[ \lambda \left(
h_{x}^{2}+h_{t}^{2}\right) +\lambda ^{3}h^{2}\right] \psi _{\lambda }dxdt+%
\frac{\gamma }{2}\left\Vert h\right\Vert _{H^{4}\left( \Omega \right) }^{2}
\\
& +B\int \displaylimits_{0}^{a}\left[ \lambda h_{x}^{2}+\lambda ^{3}h^{2}%
\right] \left( x,0\right) e^{-2\lambda x}dx.
\end{split}
\label{6.26}
\end{equation}%
Next, by (\ref{3.0}), (\ref{4.2}) and (\ref{5.2}) $R_{a,\alpha ,\mu }\subset
R_{a}\subset \Omega $ and also by (\ref{4.3}) $\psi _{\lambda }\left(
x,t\right) \geq e^{-2\lambda \left( 2\alpha a-\mu \right) }$ for $\left(
x,t\right) \in R_{a,\alpha ,\mu }.$ Hence, (\ref{6.26}) implies that%
\begin{equation*}
\begin{split}
\frac{1}{2}\int\displaylimits_{\Omega }(h_{xx}-2h_{xt})^2\psi_{\lambda}dxdt-B\int\displaylimits_{\Omega} ( h_{x}^{2}( x,0)+h^{2}) \psi _{\lambda }dxdt+\gamma \left\Vert h\right\Vert
_{H^{4}\left( \Omega \right) }^{2}  \label{6.141} \\
\geq Be^{-2\lambda \left( 2\alpha a-\mu \right) }\left( \left\Vert
h\right\Vert _{H^{1}\left( R_{a,\alpha ,\mu }\right) }^{2}+\left\Vert
h\left( x,0\right) \right\Vert _{H^{1}\left( 0,2\alpha a-\mu \right)
}^{2}\right) +\frac{\gamma }{2}\left\Vert h\right\Vert _{H^{4}\left( \Omega
\right) }^{2}. 
\end{split}
\end{equation*}
Estimates (\ref{6.210}) and (\ref{6.141}) imply the target estimate of this theorem. $\square $

\subsection{Proof of Theorem 5}

\label{sec:6.3}

Given (\ref{5.4})-(\ref{5.10}), consider $\widetilde{K}_{\lambda ,\gamma
}\left( Y^{\ast }\right) ,$ 
\begin{equation}
\widetilde{K}_{\lambda ,\gamma }\left( Y^{\ast }\right) =K_{\lambda ,\gamma
}\left( Y^{\ast }+Q\right) ,\quad Y^{\ast }\in P_{0}\left( 2d\right) .
\label{6.27}
\end{equation}%
By (\ref{4.5}), (\ref{5.301}) and (\ref{5.4}) 
\begin{equation}
M\left( Y^{\ast }+Q\right) =M\left( q^{\ast }\right) +\widetilde{M}\left(
Y^{\ast },Q-Q^{\ast }\right) =\widetilde{M}\left( Y^{\ast },Q-Q^{\ast
}\right).
\label{700}
\end{equation}%
Since functions $Y^{\ast },Q,Q^{\ast }\in H^{4}\left( \Omega \right) $ and
since by embedding theorem $H^{4}\left( \Omega \right) \subset C^{2}\left( 
\overline{\Omega }\right) ,$ then (\ref{5.3}) implies that $\left\vert 
\widetilde{M}\left( Y^{\ast },Q-Q^{\ast }\right) \right\vert \left(
x,t\right) \leq B\sigma .$ Hence, using (\ref{4.7}), (\ref{6.27}) and (\ref{700}), we obtain%
\begin{equation}
\widetilde{K}_{\lambda ,\gamma }\left( Y^{\ast }\right) \leq B\left( \sigma
^{2}+\gamma \right) .  \label{6.28}
\end{equation}%
By Theorem 4, for $Y\in P_{0}\left( 2d\right) ,$ we consider the number $%
\lambda _{2}=\lambda _{1}\left( \alpha ,3d,\Omega \right) \geq \lambda
_{1}\left( \alpha ,d,\Omega \right) .$ Recall that numbers $\lambda =\lambda
\left( \sigma \right) \geq \lambda _{2}$ and $\gamma =\gamma \left( \sigma
\right) =2e^{-\lambda \alpha \widetilde{T}}$ are the same as in (\ref{5.13})
and (\ref{5.14}) respectively. Using (\ref{500}), we obtain 
\begin{equation}
\begin{split}
& \widetilde{K}_{\lambda ,\gamma }\left( Y^{\ast }\right) -\widetilde{K}%
_{\lambda ,\gamma }\left( Y_{\min ,\lambda ,\gamma }\right) -\widetilde{K}%
_{\lambda ,\gamma }^{\prime }\left( Y_{\min ,\lambda ,\gamma }\right) \left(
Y^{\ast }-Y_{\min ,\lambda ,\gamma }\right)  \\
\geq & Be^{-2\lambda \left( 2\alpha a-\mu \right) }\left( \left\Vert Y^{\ast
}-Y_{\min ,\lambda ,\gamma }\right\Vert _{H^{1}\left( R_{a,\alpha ,\mu
}\right) }^{2}+\left\Vert Y^{\ast }\left( x,0\right) -Y_{\min ,\lambda
,\gamma }\left( x,0\right) \right\Vert _{H_{\left( 0,2\alpha a-\mu \right)
}^{1}}^{2}\right)  \\
& +\frac{\gamma }{2}\left\Vert Y^{\ast }-Y\right\Vert _{H^{4}\left( \Omega
\right) }^{2}.
\end{split}
\label{6.29}
\end{equation}%
Since $-\widetilde{K}_{\lambda ,\gamma }\left( Y_{\min ,\lambda ,\gamma
}\right) \leq 0$ and since by (\ref{5.11}) $-\widetilde{K}_{\lambda ,\gamma
}^{\prime }\left( Y_{\min ,\lambda ,\gamma }\right) \left( Y^{\ast }-Y_{\min
,\lambda ,\gamma }\right) \leq 0,$ then, \newline
using (\ref{6.28}), we estimate the first line of (\ref{6.29}) from the
above as:%
\begin{equation*}
\widetilde{K}_{\lambda ,\gamma }\left( Y^{\ast }\right) -\widetilde{K}%
_{\lambda ,\gamma }\left( Y_{\min ,\lambda ,\gamma }\right) -\widetilde{K}%
_{\lambda ,\gamma }^{\prime }\left( Y_{\min ,\lambda ,\gamma }\right) \left(
Y^{\ast }-Y_{\min ,\lambda ,\gamma }\right) \leq B\left( \sigma
^{2}+e^{-\lambda \alpha \widetilde{T}}\right) .  \label{6.30}
\end{equation*}%
We now estimate from the below the second line of (\ref{6.29}). It follows
from (\ref{5.3}), (\ref{5.4}), (\ref{5.15}) and the triangle inequality that%
\begin{equation}
\begin{split}
\left\Vert Y^{\ast }-Y_{\min ,\lambda ,\gamma }\right\Vert _{H^{1}\left(
R_{a,\alpha ,\mu }\right) }& =\left\Vert \left( Y^{\ast }+Q^{\ast }\right)
-\left( Y_{\min ,\lambda ,\gamma }+Q\right) +\left( Q-Q^{\ast }\right)
\right\Vert _{H^{1}\left( R_{a,\alpha ,\mu }\right) } \\
& \geq \left\Vert q^{\ast }-q_{\min ,\lambda ,\gamma }\right\Vert
_{H^{1}\left( R_{a,\alpha ,\mu }\right) }-\sigma .
\end{split}
\label{6.300}
\end{equation}%
It follows from the Young's inequality that $\left( x_{1}-x_{2}\right)
^{2}\geq x_{1}^{2}/2-x_{2}^{2}$ for all $x_{1},x_{2}\in \mathbb{R}.$ Hence, (%
\ref{6.300}) implies that 
\begin{equation*}
\left\Vert Y^{\ast }-Y_{\min ,\lambda ,\gamma }\right\Vert _{H^{1}\left(
R_{a,\alpha ,\mu }\right) }^{2}\geq \frac{1}{2}\left\Vert q^{\ast }-q_{\min
,\lambda ,\gamma }\right\Vert _{H^{1}\left( R_{a,\alpha ,\mu }\right)
}^{2}-\sigma ^{2}.  \label{6.31}
\end{equation*}%
Similarly, using (\ref{5.3}), (\ref{5.16}) and the trace theorem, we obtain%
\begin{equation}
\left\Vert Y^{\ast }\left( x,0\right) -Y_{\min ,\lambda ,\gamma }\left(
x,0\right) \right\Vert _{H_{\left( 0,2\alpha a-\mu \right) }^{1}}^{2}\geq
B\left\Vert r^{\ast }\left( x\right) -r_{\min ,\lambda ,\gamma }\left(
x\right) \right\Vert _{L_{2\left( 0,2\alpha a-\mu \right) }}^{2}-B\sigma
^{2}.  \label{6.32}
\end{equation}%
Combining (\ref{6.29})-(\ref{6.32}), we obtain%
\begin{equation*}
\begin{split}
& e^{-2\lambda \left( 2\alpha a-\mu \right) }\bigg(\left\Vert q^{\ast
}-q_{\min ,\lambda ,\gamma }\right\Vert _{H^{1}\left( R_{a,\alpha ,\mu
}\right) }^{2}+\left\Vert r^{\ast }\left( x\right) -r_{\min ,\lambda ,\gamma
}\left( x\right) \right\Vert _{L_{2}\left( 0,2\alpha a-\mu \right) }^{2}%
\bigg) \\
& \leq B(\sigma ^{2}+e^{-\lambda \alpha \widetilde{T}}).
\end{split}%
\end{equation*}%
Or, equivalently, 
\begin{equation}
\begin{split}
& \left\Vert q^{\ast }-q_{\min ,\lambda ,\gamma }\right\Vert _{H^{1}\left(
R_{a,\alpha ,\mu }\right) }^{2}+\left\Vert r^{\ast }\left( x\right) -r_{\min
,\lambda ,\gamma }\left( x\right) \right\Vert _{H_{\left( 0,2\alpha a-\mu
\right) }^{1}}^{2} \\
& \leq B\sigma ^{2}e^{2\lambda \left( 2\alpha a-\mu \right) }+B\exp
[-\lambda (\alpha (\widetilde{T}-4a)+\mu )].
\end{split}
\label{6.33}
\end{equation}%
Recall that $\ln \sigma _{0}^{-1/\left( 2\left( 2\alpha a-\mu \right)
\right) }\geq \lambda _{2}$ and $\sigma \in \left( 0,\sigma _{0}\right) .$
Since by (\ref{5.13}) $\lambda =\lambda \left( \sigma \right) >\lambda _{2},$
then in (\ref{6.33}) 
\begin{equation}
\sigma ^{2}e^{2\lambda \left( 2\alpha a-\mu \right) }=\sigma \text{, }\exp %
\left[ -\lambda \left( \alpha \left( \widetilde{T}-4a\right) +2\mu \right) %
\right] =\sigma ^{\nu },  \label{6.34}
\end{equation}%
where $\nu >0$ is the number defined in (\ref{5.12}). Thus, using (\ref{6.33}%
) and (\ref{6.34}), we obtain%
\begin{align}
& \left\Vert q^{\ast }-q_{\min ,\lambda ,\gamma }\right\Vert _{H^{1}\left(
R_{b,\alpha ,\mu }\right) }\leq B\left( \sqrt{\sigma }+\sigma ^{\nu
/2}\right) ,  \label{6.35} \\
& \left\Vert r^{\ast }\left( x\right) -r_{\min ,\lambda ,\gamma }\left(
x\right) \right\Vert _{L_{2\left( 0,2\alpha b-\mu \right) }}\leq B\left( 
\sqrt{\sigma }+\sigma ^{\nu /2}\right) .  \label{6.36}
\end{align}%
Estimates (\ref{6.35}) and (\ref{6.36}) combined with (\ref{5.12}) imply the
target estimates (\ref{5.17}) and (\ref{5.19}) of this theorem. \ \ \ $%
\square $

\subsection{Proof of Theorem 6}

\label{sec:6.4}

Combining Theorem 2 with Theorem 2.1 of \cite{KlibanovNik:ra2017}, we obtain
that the number $\theta \in \left( 0,1\right) $ exists and estimate (\ref%
{5.22}) holds. Estimate (\ref{5.23}) follows from (\ref{5.22}), the trace
theorem and (\ref{3.6}). Estimate (\ref{5.24}) follows from (\ref{5.17}), (%
\ref{5.22}) and the triangle inequality. Similarly, (\ref{5.25}) follows
from (\ref{5.19}), (\ref{5.23}) and the triangle inequality. $\ \square $

\section{Numerical Implementation}

\label{sec:7}

In this section, we describe our numerical procedure to solve the
Minimization Problem formulated in section 4.2. We work with the finite
difference analog of the functional $K_{\lambda ,\gamma }$ defined in (\ref%
{4.7}). To do this, we use the uniform grid $\Omega _{h}\subset \Omega ,$ where $h=(h_{x},h_{t})$. More
precisely, for certain integers $N_{x},N_{t}>1$ 
\begin{equation}
\Omega _{h}=\{(x,t):x=(i-1)h_{x},\hspace{0.3em}t=(j-1)h_{t},\hspace{0.3em}i=1,\dots
,N_{x}+1,\hspace{0.3em}j=1,\dots ,N_{t}+1\}.
\label{7.1}
\end{equation}
In all numerical studies of this paper we take $a=1.1\sqrt{\overline{c}},$ see (\ref{2.16}). Therefore, by (\ref{3.0}) we have $\Omega $ as
\begin{equation*}
\Omega =\left( 0,1.1\sqrt{\overline{c}}\right) \times (0,2.2\sqrt{\overline{c}}).
\end{equation*}

We denote values of the function functions $q\left( x,t\right) $ at the grid
points of the domain $\Omega _{h}$ by $q_{i,j}=q(x_{i},t_{j})$. Define the
finite-difference analog $M_{i,j}=M(q_{i,j})$ of the operator $M(q)$ in (\ref%
{4.5}) at the point $(x_{i},t_{j})$ as 
\begin{equation}
\begin{split}
M_{i,j}& =\left( \frac{q_{i-1,j}-2q_{i,j}+q_{i+1,j}}{h_{x}^{2}}\right)
-2\left( \left( \frac{q_{i+1,j+1}-q_{i+1,j}}{h_{x}h_{t}}\right) -\left( 
\frac{q_{i,j+1}-q_{i,j}}{h_{x}h_{t}}\right) \right) \\
& +4\left( \frac{q_{i+1,1}-q_{i,1}}{h_{x}}\right) q_{i,j},
\end{split}%
\end{equation}%
Following (\ref{3.6}), denote 
\begin{equation}
r_{comp}(x)=4\left( \frac{q_{i+1,1}-q_{i,1}}{h_{x}}\right) ,\quad \text{ for 
}x=x_{i},\hspace{0.3em}\forall i\in \lbrack 1,N_{x}].  \label{r_comp}
\end{equation}%
Thus, we define the finite difference analog $K_{\lambda ,\gamma }^{h}\left(
q\right) $ of the functional $K_{\lambda ,\gamma }$ as 
\begin{equation}
\begin{split}
&K_{\lambda ,\gamma }^{h}\left( q\right)
=\sum_{i=2}^{N_{x}}\sum_{m=1}^{N_{t}}\left[ M_{i,j}\right] ^{2}\psi
_{\lambda }(x_{i},t_{j})h_{x}h_{t}+\gamma
\sum_{i=1}^{N_{x}+1}\sum_{m=1}^{N_{t}+1}(q_{i,j})^{2} h_x h_t \\
& +\gamma \sum_{i=1}^{N_{x}}\sum_{m=1}^{N_{t}}\Bigg\{\left( \frac{q_{i+1,j}-q_{i,j}}{h_{x}}\right) ^{2}+\left( \frac{q_{i,j+1}-q_{i,j}}{h_{t}}\right) ^{2} \Bigg\} h_x h_t \\ &+\gamma \sum_{i=2}^{N_{x}}\sum_{m=2}^{N_{t}}\Bigg\{ \left( \frac{q_{i-1,j}-2q_{i,j}+q_{i+1,j}}{h_{x}^{2}}\right) ^{2}+\left( \frac{q_{i,j-1}-2q_{i,j}+q_{i,j+1}}{h_{t}^{2}}\right) ^{2}\Bigg\} h_{x}h_{t},
\end{split}
\label{7.4}
\end{equation}
Supplying the gradient of the functional $K_{\lambda ,\gamma }^{h}$ via explicit formula significantly reduces computational time of the minimization procedure. The gradient of the functional $\nabla K_{\lambda ,\gamma }^{h}$ is defined as a vector of partial derivatives with respect to the values $%
q_{i,j}$ of the function at the grid points $x=x_{i^{\ast }},t=t_{j^{\ast }}$%
\begin{equation}
\begin{split}
&\frac{\partial K_{\lambda ,\gamma }^{h}\left( q\right) }{\partial
q_{i^{\ast },j^{\ast }}}= \sum_{i=2}^{N_{x}}\sum_{m=1}^{N_{t}} 2 M_{i,j} \frac{\partial M_{i,j}}{\partial
q_{i^{\ast },j^{\ast }}} \psi
_{\lambda }(x_{i},t_{j})h_{x}h_{t}
+2 \gamma \Bigg\{\sum_{i=1}^{N_{x}+1}\sum_{m=1}^{N_{t}+1} q_{i,j} \frac{\partial q_{i,j}}{\partial
q_{i^{\ast },j^{\ast }}}\\ &+\sum_{i=1}^{N_{x}}\sum_{m=1}^{N_{t}}\frac{1}{h_x}\left( \frac{%
q_{i+1,j}-q_{i,j}}{h_{x}}\right) \left(\frac{\partial q_{i+1,j}}{\partial
q_{i^{\ast },j^{\ast }}}-\frac{\partial q_{i,j}}{\partial
q_{i^{\ast },j^{\ast }}}\right) \\
&+\sum_{i=1}^{N_{x}}\sum_{m=1}^{N_{t}}\frac{1}{h_t}\left( \frac{q_{i,j+1}-q_{i,j}}{h_{t}}\right)\left(\frac{\partial q_{i,j+1}}{\partial
q_{i^{\ast },j^{\ast }}}-\frac{\partial q_{i,j}}{\partial
q_{i^{\ast },j^{\ast }}}\right)\\
&+\sum_{i=2}^{N_{x}}\sum_{m=2}^{N_{t}}\frac{1}{h_x^2}\left( \frac{%
q_{i-1,j}-2q_{i,j}+q_{i+1,j}}{h_{x}^{2}}\right)\left(\frac{\partial q_{i-1,j}}{\partial
q_{i^{\ast },j^{\ast }}}-2\frac{\partial q_{i,j}}{\partial
q_{i^{\ast },j^{\ast }}}+\frac{\partial q_{i+1,j}}{\partial
q_{i^{\ast },j^{\ast }}}\right)\\
&+\sum_{i=2}^{N_{x}}\sum_{m=2}^{N_{t}}\frac{1}{h_t^2}\left( \frac{%
q_{i,j-1}-2q_{i,j}+q_{i,j+1}}{h_{t}^{2}}\right)\left(\frac{\partial q_{i,j-1}}{\partial
q_{i^{\ast },j^{\ast }}}-2\frac{\partial q_{i,j}}{\partial
q_{i^{\ast },j^{\ast }}}+\frac{\partial q_{i,j+1}}{\partial
q_{i^{\ast },j^{\ast }}}\right) \Bigg\}h_{x}h_{t}.
\end{split}
\label{7.50}
\end{equation}
However the formula (\ref{7.50}) is not an explicit one yet.
Hence, we use the formula $\frac{\partial q_{i,j}}{\partial
q_{i^{\ast },j^{\ast }}} = \delta_{i^{\ast},j^{\ast}}$, where $\delta_{i,j}$ denotes Kronecker symbol. 
\newline

The following simplest example explains our idea of using Kronecker symbol on how to obtain explicit formulas of derivatives in (\ref{7.50})  
\begin{equation*}
\begin{split}
&\frac{\partial}{\partial
q_{i^{\ast },j^{\ast }}} \sum_{i=1}^{N_{x}}\sum_{m=1}^{N_{t}} \left(\frac{%
q_{i+1,j}-q_{i,j}}{h_{x}} \right)^2 = \sum_{i=1}^{N_{x}}\sum_{m=1}^{N_{t}} \frac{2}{h_x} \left(\frac{q_{i+1,j}-q_{i,j}}{h_{x}} \right) \left(\frac{\partial q_{i+1,j}}{\partial
q_{i^{\ast },j^{\ast }}}-\frac{\partial q_{i,j}}{\partial
q_{i^{\ast },j^{\ast }}}\right) = \\
&\sum_{i=1}^{N_{x}}\sum_{m=1}^{N_{t}} \frac{2}{h_x} \left(\frac{q_{i+1,j}-q_{i,j}}{h_{x}} \right) \left(\delta_{i^{\ast}-1,j^{\ast}}-\delta_{i^{\ast},j^{\ast}}\right) = -\frac{2}{h_x} \left(\frac{q_{i^{\ast}-1,j^{\ast}}-2q_{i^{\ast},j^{\ast}}+q_{i^{\ast}+1,j^{\ast}}}{h_{x}} \right). 
\end{split}
\end{equation*}
Thus we obtain the following explicit formulas for the derivatives in (\ref{7.50})
\begin{equation}
\begin{split}
& \frac{\partial K_{\lambda ,\gamma }^{h}\left( q\right) }{\partial
q_{i^{\ast },j^{\ast }}}=h_{x}h_{t}\bigg\{\frac{2}{h_{x}^{2}}\left(
M_{i^{\ast }-1,j^{\ast }}\psi _{\lambda }^{i^{\ast }-1,j^{\ast
}}-2M_{i^{\ast },j^{\ast }}\psi _{\lambda }^{i^{\ast },j^{\ast }}+M_{i^{\ast
}+1,j^{\ast }}\psi _{\lambda }^{i^{\ast }+1,j^{\ast }}\right) \\
& -\frac{4}{h_{x}h_{t}}\left( M_{i^{\ast }-1,j^{\ast }-1}\psi _{\lambda
}^{i^{\ast }-1,j^{\ast }-1}-M_{i^{\ast }-1,j^{\ast }}\psi _{\lambda
}^{i^{\ast }-1,j^{\ast }}-M_{i^{\ast },j^{\ast }-1}\psi _{\lambda }^{i^{\ast
},j^{\ast }-1}+M_{i^{\ast },j^{\ast }}\psi _{\lambda }^{i^{\ast },j^{\ast
}}\right) \\
& +\frac{8}{h_{x}}\left( M_{i^{\ast }-1,1}\psi _{\lambda }^{i^{\ast
}-1,1}q_{i^{\ast }-1,1}-M_{i^{\ast },1}\psi
_{\lambda }^{i^{\ast },1}q_{i^{\ast },1}+M_{i^{\ast
},j^{\ast }}\psi _{\lambda }^{i^{\ast },j^{\ast }}(q_{i^{\ast
}+1,1}-q_{i^{\ast },1})\right) \bigg\} \\
& +2\gamma h_{x}h_{t}\bigg\{\frac{q_{i^{\ast }+2,j^{\ast }}-4q_{i^{\ast
}+1,j^{\ast }}+6q_{i^{\ast },j^{\ast }}-4q_{i^{\ast }-1,j^{\ast
}}+q_{i^{\ast }-2,j^{\ast }}}{h_{x}^{4}} \\
&- \frac{q_{i^{\ast }+1,j^{\ast }}-2q_{i^{\ast },j^{\ast }}+q_{i^{\ast
}-1,j^{\ast}}}{h_{x}^{2}} +\frac{q_{i^{\ast
},j^{\ast }+2}-4q_{i^{\ast },j^{\ast }+1}+6q_{i^{\ast },j^{\ast
}}-4q_{i^{\ast },j^{\ast }-1}+q_{i^{\ast },j^{\ast }-2}}{h_{t}^{4}}\\
&-\frac{q_{i^{\ast },j^{\ast}+1}-2q_{i^{\ast },j^{\ast}}+q_{i^{\ast },j^{\ast}-1}}{h_{t}^{2}}
+q_{i^{\ast },j^{\ast }}\bigg\},
\end{split}
\label{7.60}
\end{equation}%
where $i^{\ast }\in \lbrack 3,N_{x}-1],\hspace{0.3em}j^{\ast }\in \lbrack
3,N_{t}-1]$ and $\psi _{\lambda }^{i^{\ast },j^{\ast }}=\psi _{\lambda
}(x_{i^{\ast }},t_{j^{\ast }})$.

The derivatives in (\ref{7.60}) are defined only in the interior of the domain $\Omega^h$, defined in (\ref{7.1}). In addition it is necessary to approximate the partial derivatives of $K_{\lambda,\gamma }^{h}\left( q\right) $ on the boundary of this domain, using Taylor
expansion of $\partial K_{\lambda ,\gamma }^{h}(q)/\partial q_{i^{\ast },j}$. This is because the function $r_{comp}(x)$
defined in (\ref{r_comp}) is updated on every iteration of the gradient
descent method. Thus, for $i^{\ast }=1,2$ and $j\in \lbrack
3,N_{x}-1]$ we set
\begin{equation}
\begin{split}
&\frac{\partial K^h_{\lambda ,\gamma }\left( q\right)}{\partial
q_{i^{\ast},j^{\ast}}} = \frac{\partial K^h_{\lambda ,\gamma }\left( q\right)}{\partial
q_{i^{\ast}+1,j^{\ast}}} - \left(  \frac{1}{h_x} \left(\frac{\partial K^h_{\lambda ,\gamma }\left( q\right)}{\partial
q_{i^{\ast}+2,j^{\ast}}} - \frac{\partial K^h_{\lambda ,\gamma }\left( q\right)}{\partial
q_{i^{\ast}+1,j^{\ast}}}\right)\right) h_x + \Bigg( \frac{1}{h_x^2}\Bigg(\frac{\partial K^h_{\lambda ,\gamma }\left( q\right)}{\partial
q_{i^{\ast}+1,j^{\ast}}}\\ 
&- 2\frac{\partial K^h_{\lambda ,\gamma }\left( q\right)}{\partial
q_{i^{\ast}+1,j^{\ast}}} 
+\frac{\partial K^h_{\lambda ,\gamma }\left( q\right)}{\partial
q_{i^{\ast}+3,j^{\ast}}}\Bigg) \Bigg) h_x^2 + o(h_x^2) \approx \frac{5}{2} \frac{\partial K^h_{\lambda ,\gamma
}\left( q\right)}{\partial q_{i^{\ast}+1,j}} - 2 \frac{\partial K^h_{\lambda
,\gamma }\left( q\right)}{\partial q_{i^{\ast}+2,j}} + \frac{1}{2} \frac{%
\partial K^h_{\lambda ,\gamma }\left( q\right)}{\partial q_{i^{\ast}+3,j}}.
\label{7.3}
\end{split}
\end{equation}
Thereafter for $i^{\ast} = N_x,N_x+1$ and $j \in [3,N_x-1]$
\begin{equation*}
\begin{split}
\frac{\partial K^h_{\lambda ,\gamma }\left( q\right)}{\partial
q_{i^{\ast},j^{\ast}}} = \frac{5}{2} \frac{\partial K^h_{\lambda ,\gamma
}\left( q\right)}{\partial q_{i^{\ast}-1,j}} - 2 \frac{\partial K^h_{\lambda
,\gamma }\left( q\right)}{\partial q_{i^{\ast}-2,j}} + \frac{1}{2} \frac{%
\partial K^h_{\lambda ,\gamma }\left( q\right)}{\partial q_{i^{\ast}-3,j}} ,
\label{7.8}
\end{split}
\end{equation*}
we define $\partial K_{\lambda ,\gamma }^{h}\left( q\right) /\partial
q_{i,j^{\ast }}$ on the rest of the boundary, i.e. for $j^{\ast
}=1,2,N_{t},N_{t}+1$ and \newline
$i\in \lbrack 1,N_{x}+1]$, similarly to (\ref{7.4}),(\ref{7.3}).

\begin{remark}
Note that the functional in (\ref{7.4}) uses the Tikhonov regularization term in the $H^2(\Omega_h)$ norm instead of the $H^3(\Omega_h)$ required by the theory. We have established numerically that this is sufficient for computations.   
\end{remark}

\subsection{Gradient descent method (GDM)}

\label{sec:7.1}

Even though Theorem 6 guarantees global convergence of the gradient projection method, we have numerically established that the simpler to implement gradient descent method (GDM) works well for our studies. The latter coincides with the conclusions of all above cited publications on the numerical studies of the convexification, e.g. see \cite{KlibanovNik:ra2017,klibanov2019backscatter,klibanovzhang2020convexification,klibanov2019convexification,klibanov2018new,klibanov2017globally}. We now apply the GDM to find the minimizer of functional (\ref{7.4}). According to our theory, the initial guess for the GDM can be an arbitrary function $q^{0}(x,t)\in P(d,s_{0},s_{1})$.
We choose the initial guess $q^{0}(x,t)$ as the solution to the following
problem, derived from (\ref{3.4})-(\ref{3.5}) by setting $r(x)=4q_{x}(x,0)=0$
and assuming $s_{1}(t)=0$ for $t\geq 2b,$ see (\ref{3.0}). Thus,
\begin{align*}
& q_{xx}^{0}-2q_{xt}^{0}=0,\quad \left( x,t\right) \in (0,b)\times (0,%
\widetilde{T}),  \\
& q^{0}\left( 0,t\right) =s_{0}\left( t\right) ,\quad q_{x}^{0}\left(
0,t\right) =s_{1}\left( t\right) ,\quad q_{x}^{0}(b,t)=0,\quad t\in (0,%
\widetilde{T}),
\end{align*}
which has the unique solution 
\begin{equation}
q^{0}(x,t)=s_{0}(t)+\frac{1}{2}\int \displaylimits_{t}^{t+2x}s_{1}(\tau
)d\tau , \quad (x,t) \in \Omega,  \label{q0}
\end{equation}%
where functions $s_{0}(t)$ and $s_{1}(t)$ are obtained from the pre-processed data $g_{0}(t),g_{1}(t)$ (\ref{2.5}), see Appendix for a more detailed explanation of the data simulation and pre-processing. In the case of the computationally simulated data
the functions $g_{0}(t),g_{1}(t)$ are obtained from the numerical
solution of the Forward Problem (\ref{2.3})-(\ref{2.4}) via the tridiagonal matrix algorithm \cite{conte2017elementary}.

\begin{remark}
Other forms of the initial guess in a discrete analog $P^h(d,s_{0},s_{1})$ of the set $P(d,s_{0},s_{1})$ can also be chosen, see the section below. Such choices  will not lead to significant different solutions, due to the global convergence of our numerical method, see Theorem 6 of section \ref{sec:5}.  
\end{remark}
In all numerical tests of this paper we choose $N_{x}=100,N_{t}=100$. In all further computations of inverse problems with simulated data we use the multiplicative random noise of the level $\delta ^{\ast }=0.05,$ i.e. $5\%$. Since we use functions $s_{0}(t)=g_{0}^{\prime }(t)$ and $s_{1}(t)=g_{0}^{\prime \prime }(t)+g_{1}^{\prime }(t),$ we need to differentiate the noisy data $g_{0}(t),g_{1}(t).$, which we do by taking derivatives of the envelopes, see Figure \ref{fig2}. More detailed decription of data pre-processing and differentiation is given in Appendix.
\begin{figure}[htb!]
\begin{center}
\subfloat[$g_0(t)$ and its envelope.]{\includegraphics[width
=.45\textwidth]{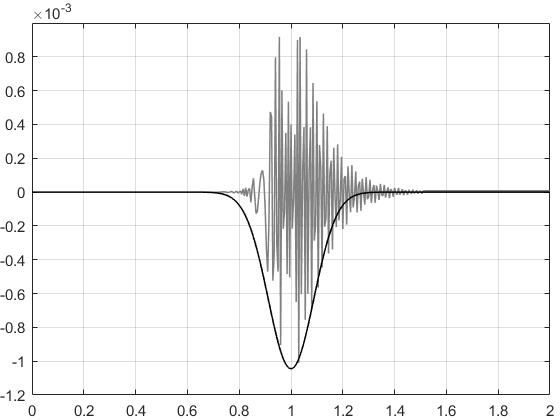}} \quad 
\subfloat[$g_1(t)$ and its envelope.]{\includegraphics[width =.45\textwidth]{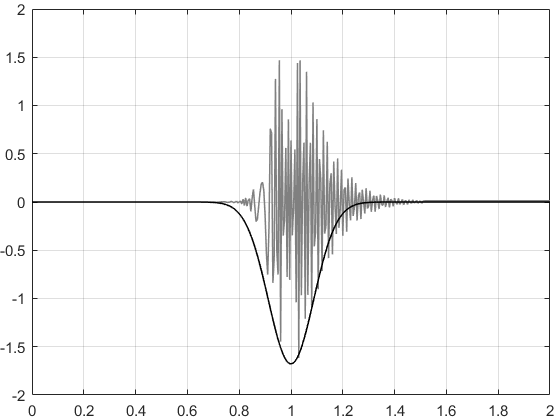}}
\end{center}
\caption{\emph{The simulated data $g_0(t),g_1(t)$ for $\widetilde{c}(y)$ in (\ref{7.10}), $w^{\ast} = 0.1$, $A_c = 0.2$ and the multiplicative random noise
level $\delta =0.05$. The solid line depicts the data computed from the solution of the Forward Problem with $5\%$ added multiplicative random noise. }}
\label{fig2}
\end{figure} \newline
The following values of parameters were used:%
\begin{equation}
\lambda =2,\quad \gamma =10^{-6},\quad \alpha =0.5, \quad h_x = 0.01, \quad h_t = 0.02.  \label{7.100}
\end{equation}
The parameters (\ref{7.100}) were found by the trial and error procedure, and they work well for the numerical studies of this paper. We point out that even though our above theory requires large values of the parameter $\lambda ,$ our numerical experience tells us that $\lambda =2$ provides decent results. This is consistent with
other numerical results on the convexification, where it was computationally established that $\lambda \in \lbrack 1,3]$ works well numerically \cite{KlibanovNik:ra2017,klibanov2019backscatter,klibanovzhang2020convexification,klibanov2019convexification,klibanov2018new,klibanov2017globally}. The topic of optimal choices of parameters is outside of the scope of this paper. For brevity we use everywhere below notations for functions and variables in continuous setting, although we work with discrete setting as stated above.

Denote the result obtained on $k$-th iteration of the GDM by $%
q^{k}(x,t)$. Using this function, we calculate the function $r^{k}(x)$
via (\ref{3.6}). After a sufficient number of iterations $k_{stop}$ we come up with the computed coefficient $r_{comp}(x)=r^{k_{stop}}(x)$. Thus we formulate Algorithm \ref{alg 1} for solving the Minimization Problem for functional $K_{\lambda ,\gamma
}^{h}(q)$ on the set $P^{h}(d,s_{0},s_{1})$,  i.e. in (\ref{4.6}) $%
\Omega $ is replaced with $\Omega _{h},\ H^{3}\left( \Omega \right) $ is
replaced with the discrete analog of $H^{2}\left( \Omega \right).$ We use the discrete version of the space $H^{2}\left( \Omega \right) $ because the regularization term in $K_{\lambda ,\gamma }^{h}\left( q\right) $ uses the discrete form of the norm in this space.

\begin{algorithm}
\caption{The minimization of the functional $K^h_{\lambda ,\gamma }$ via GDM.}
\begin{algorithmic}[1]
    \State\, Compute function $q^0(x,t)$ via (\ref{q0}).
    \State\, Compute $\nabla K^h_{\lambda ,\gamma }(q^k)$  and perform one step of the GDM to find $q^{k+1}(x,t), k = 0,1, \dots$.
     \State\, Calculate $r_{k+1}(x)$ using $q^{k+1}(x,t)$ in (\ref{r_comp}).
    \State\, Repeat steps 2,3 while $\Vert K^h_{\lambda,\gamma} (q^k)\Vert_{\infty} > 10^{-2}\Vert K^h_{\lambda,\gamma} (q^0)\Vert_{\infty}$ or \newline $\Vert \nabla K^h_{\lambda,\gamma} (q^k) \Vert_{\infty} > 10^{-2}\Vert \nabla K^h_{\lambda,\gamma} (q^0)\Vert_{\infty}$. Let $k_{stop}$ be the last iteration number at which either of these inequalities holds. Then stop at $k=k_{stop}.$.
    \State\, Compute the $r_{comp}(x)$ by  (\ref{r_comp}) using $q^{k_{stop}}(x,t)$.
\end{algorithmic}
\label{alg 1}
\end{algorithm}

\begin{figure}[htb!]
\begin{center}
\subfloat[Test 1. True $\widetilde{c}(y)$ given by (\ref{7.10}) with $A_c = 0.2$, \newline $w^{\ast} = 0.075.$ The horizontal axis depicts the values of $y$.]{\includegraphics[width
=.45\textwidth]{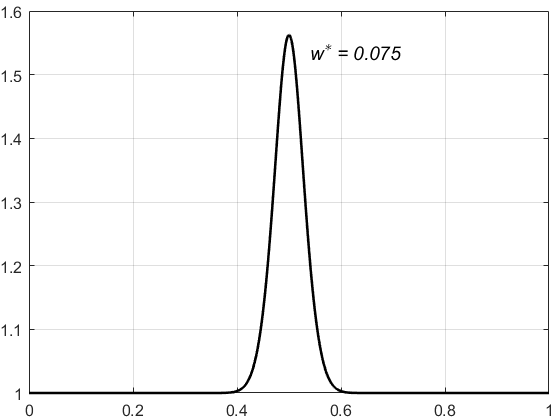}} \quad 
\subfloat[Test 1. $r_{comp}(x)$ computed via the Algorithm \ref{alg 1} compared to the true function $r^{\ast}(x)$ defined by $\widetilde{c}(y)$ in (\ref{7.10}) depicted on (a). The horizontal axis depicts the values of $x$.]{\includegraphics[width
=.45\textwidth]{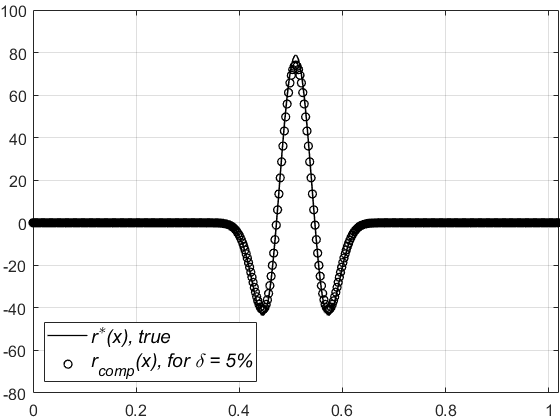}} \\
\subfloat[Test 2. True $\widetilde{c}(y)$ given by the sum of two functions of the form (\ref{7.10}), centered at $y = 0.3$ and $y= 0.7$ with $A_c = 0.2$ and $w^{\ast}_1 = 0.1, w^{\ast}_2 = 0.1$. ]{\includegraphics[width
=.45\textwidth]{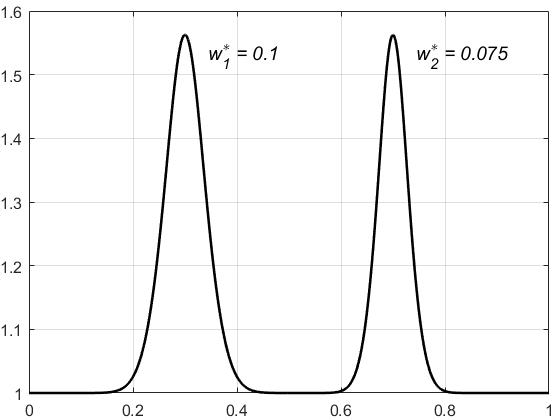}} \quad\subfloat[Test 2.  $r_{comp}(x)$ computed via the Algorithm \ref{alg 1} compared to the true function $r^{\ast}(x)$ defined by $\widetilde{c}(y)$ in (\ref{7.10}) depicted on (b). The horizontal axis depicts the values of $x$.]{\includegraphics[width
=.45\textwidth]{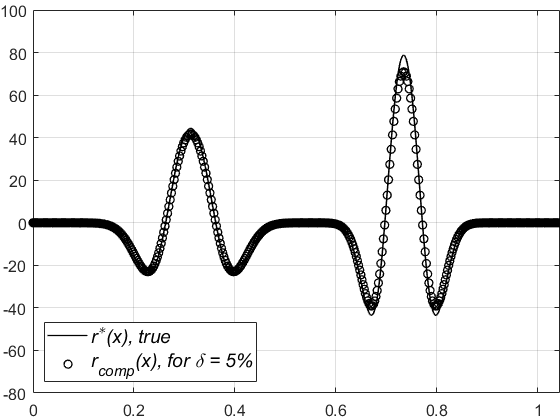}}
\end{center}
\caption{\emph{Numerical solutions (b),(d) of CIP2 for noisy data with $\delta = 0.05$ for $\widetilde{c}(y)$ depicted on (a),(c). See Table 1 for the reconstruction errors.}}
\label{fig3}
\end{figure}

Thus, the numerical procedure described in Algorithm \ref{alg 1} delivers the numerical solution of the CIP2. We describe here two tests for the developed numerical algorithm. The calculated functions $r_{comp}(x)$ for both tests were interpolated with cubic splines on the finer grid with $N'_x = 450$ grid points. This interpolation is an important step for the reconstruction of the target functions $c(y)$, see section \ref{sec:7.2}. The results of the reconstruction of the functions $r_{comp}(x)$ are depicted on Figure \ref{fig3}.
The corresponding errors for the reconstructions are summarized in Table 1. 

\textbf{Test 1.} The true function $r^{\ast}(x)$ to be reconstructed as well as the simulated data were calculated from $\widetilde{c}(y)$ given by
\begin{equation}
\widetilde{c}(y) = \left(1-A_c \exp{\left(-\frac{(y-0.5)^2}{2 w_c^2}\right)}%
\right)^{-2}, \quad w_c = 2 \sqrt{2 \ln{2}} w^{\ast} ,
\label{7.10}
\end{equation}
where $w^{\ast}$ is full width at half maximum. The maximal value of the function $\widetilde{c}(y)$ is \\ max$(\widetilde{c}%
^{\ast}(y)) = 1/(1-A_c)^2$. For the first test we choose $A_c = 0.2$ and $w^{\ast} = 0.075$, which corresponds to the maximal value of the dielectric constant of max$(c(y)) = 1.56$, see Figure \ref{fig3} (a) and solid line on Figure \ref{fig4} (a).

\textbf{Test 2.} The true function for the second test $r^{\ast}(x)$ is calculated from $\widetilde{c}(y)$ given by the sum of two functions of the form (\ref{7.10}), centered at $y = 0.3$ and $y= 0.7$ with $A_c = 0.2$ and $w_1 = 2 \sqrt{2 \ln{2}} w_1^{\ast}, \hspace{0.3em} w^{\ast}_1 = 0.1$ and $w_2 = 2 \sqrt{2 \ln{2}} w_2^{\ast}, \hspace{0.3em} w^{\ast}_2 = 0.075$, see Figure \ref{fig3} (b) and solid line on \ref{fig6} (a).
\begin{equation*}
\widetilde{c}(y) = \left(1-A_c \exp{\left(-\frac{(y-0.3)^2}{2 w_1^2}\right)}-A_c \exp{\left(-\frac{(y-0.7)^2}{2 w_2^2}\right)}%
\right)^{-2}.
\end{equation*}
\newpage
\begin{center}
\textbf{Table 1.} Relative errors for the numerical solutions of the CIP2.
\vspace{0.5em}
\begin{tabular}{l|l}
       & $\Vert r^{\ast} - r_{comp} \Vert_{L^2(0,x(b))}/\Vert r^{\ast} \Vert_{L^2(0,x(b))}$ \\ \hline
Test 1. & 0.0600 \\ \hline
Test 2. & 0.0900 \\ \hline
\end{tabular}
\end{center}

\subsection{Reconstruction of $c(y)$ from the $r_{comp}(x)$ using weighted least-squares (WLS)} In this section we describe the second stage of the reconstruction procedure. More precisely, we show how to numerically reconstruct function $\widetilde{c}(y)$ from $r_{comp}(x)$. Since $c(y) \in [1, \overline{c}]$ and $x \in (0,a)$, then the upper estimate for the value of $y$ is $\overline{y} = $ max$(1,a)$. However, we were only interested in the values of $c(y)$ on the interval $(0,1)$, therefore we enforce $\overline{y} = 1$. Assume that the the values of the function $r\left( x\right)$ are obtained via Algorithm \ref{alg 1} of section \ref{sec:7.1} on the interval $(0,a)$. Relation between $c(y)$ and $r(x)$ is given by the quasilinear
differential equation (\ref{2.8}) coupled with (\ref{2.6}). 

Denote $\widetilde{c}(y)=c(y)^{-1/2}$. Then, using the substitution $S(x) = \sqrt{\widetilde{c}(y(x))}$ we arrive at the following initial value problem for $\widetilde{p}(x) = \widetilde{c}(y(x))$:
\vspace{1em}

\emph{Given the function $
r(x), \hspace{0.3em} x \in (0,a)$, find the function $\widetilde{p}(x) \in C^2(\mathbb{R})$ and satisfying} 
\begin{align}
&\widetilde{p}^{\prime \prime }(x)\widetilde{p}(x)/2-(\widetilde{p}^{\prime}(x))^2/4 = r(x), \hspace{0.3em} x \in (0,a),  \label{7.91} \\
& \widetilde{p}(0) = 1, \hspace{0.3em} \widetilde{p}^{\prime }(0) = 0, 
\hspace{0.3em} \widetilde{p}^{\prime \prime }(0) = 0 ,  \label{7.92} \\
&y^{\prime }(x) = \widetilde{p}(x), 
\hspace{0.3em} y(0) = 0.  \label{7.93}
\end{align}

There are two significant difficulties associated with the numerical solution of problem (\ref{7.91})-(\ref{7.93}). First, the dependence $y=y(x)$ is unknown, whereas the right hand side $r(x)$ of equation (\ref{7.91}) is reconstructed by Algorithm 7.1 as the function of $x$. 
Second, the problem (\ref{7.91})-(\ref{7.93}) is overdetermined. We explain below how do we handle these two difficulties. 

\begin{figure}[htb!]
\begin{center}
\subfloat[$r(x)$ computed via (\ref{7.91}) for $\widetilde{p}(x) = \widetilde{c}(y(x))$ in (\ref{7.10}), $w^{\ast} = 0.075$, $A_c = 0.2$, the corresponding $l = b_2-b_1 = 0.073$. In addition vertical dashed lines indicate edges of the interval $(b_1,b_2)$, where $r(x) > 0$.]{\includegraphics[width =.45\textwidth]{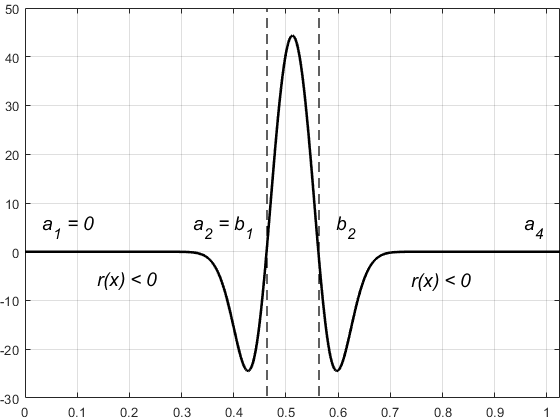}}
\quad 
\subfloat[Optimal choice of $\rho^{\ast}(l)$, computed via WLS for $\widetilde{p}(x) = \widetilde{c}(y(x))$ in (\ref{7.10}), $l \approx w^{\ast}$, $A_c = 0.2$ (triangles) and $A_c = 0.5130$ (squares) in double log scale.]{\includegraphics[width
=.45\textwidth]{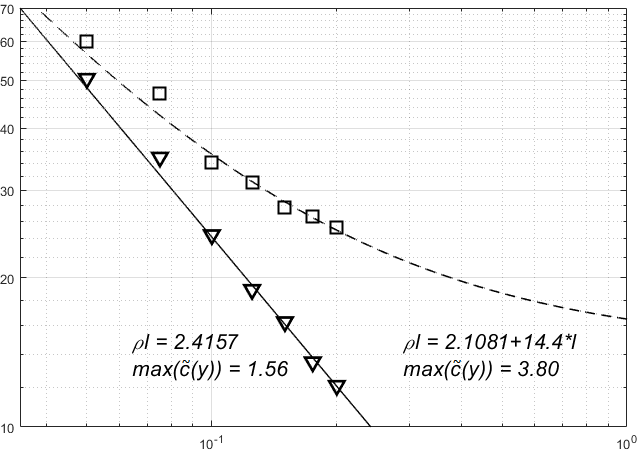}}
\end{center}
\caption{\emph{The optimal parameters of WLS reconstruction method for the
test function of the form (\protect\ref{7.10}). The horizontal for (a) depicts the values of $x$ and the horizontal axis for (b) depicts the values of $l$. }}
\label{fig4}
\end{figure}

We have attempted to solve problem (\ref{7.91})-(\ref{7.93}) via the Runge-Kutta method, on of the whole interval $x \in [0,a]$ for a piecewise constant function $\widetilde{p}(x)$ without the condition \newline $\widetilde{p}^{\prime \prime }(0) = 0$. Additionally, we assume here that $\widetilde{p}(x) = 1, \text{ for } x \leq 0$. This allows us to recover $\widetilde{p}(x)$ and to compute $y = y(x)$ via (\ref{7.93}). However, we have observed numerically that this approach does not provide sufficiently accurate reconstructions due to the fact that the function $r(x)$ attains both positive and negative values, see Figure \ref{fig4} (a).
Hence, we formulate an alternative approach to solve problem (\ref{7.91})-(\ref{7.93}). This approach requires the assumption $\widetilde{p}^{\prime \prime }(0) = 0$. Note that this assumption is quite natural, because of (\ref{2.2}).
\begin{figure}[htb!]
\begin{center}
\subfloat[$\widetilde{c}_{comp}(y)$ computed via WLS with optimal value of $\rho = 34.98$ (circles) and via WLS with $\rho= 0$ (dashed line). Solid line depicts the exact $\widetilde{p}(x)$, which is generated by $\widetilde{c}(y)$ in (\ref{7.10}).]{\includegraphics[width
=.45\textwidth]{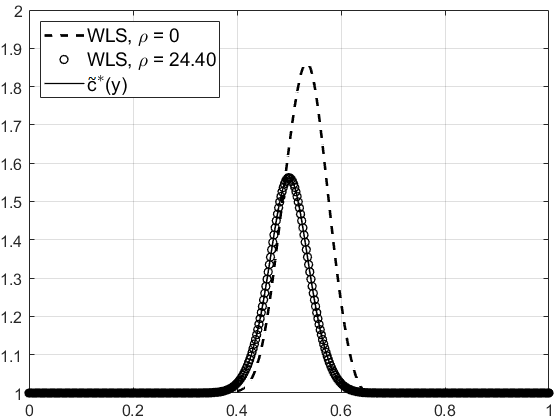}} \quad
\subfloat[$\widetilde{c}_{comp}(y)-\widetilde{c}^{\ast}(y)$ for WLS with optimal choice of $\rho = \rho^{\ast}(l)$ (crossed line) compared with the Runge-Kutta method (dashed line).]{\includegraphics[width
=.45\textwidth]{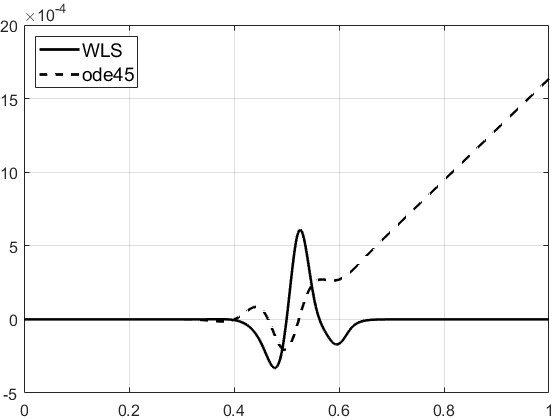}} 
\end{center}
\caption{\emph{The horizontal axis depicts the values of $y$, see the text of section \ref{sec:7.2} for an explanation. (a) Computed functions and  (b) reconstruction error $\protect\widetilde{c}_{comp}(y)-\protect%
\widetilde{c}^{\ast}(y)$ for WLS and for the MATLAB built-in "ode45", a
six-stage, fifth-order, Runge-Kutta method, where $\protect\widetilde{c}%
_{comp}(y)$ and $\protect\widetilde{c}^{\ast}(y)$ denote the computed and
true functions respectively. }}
\label{fig5}
\end{figure}

\begin{figure}[htb!]
\begin{center}
\subfloat[$\widetilde{c}_{comp}(y)$,
computed via WLS with ($\rho = 34.98$, crossed line) and without ($\rho  =
0$, solid line) an optimal choice of $\rho(l)$.]{\includegraphics[width
=.45\textwidth]{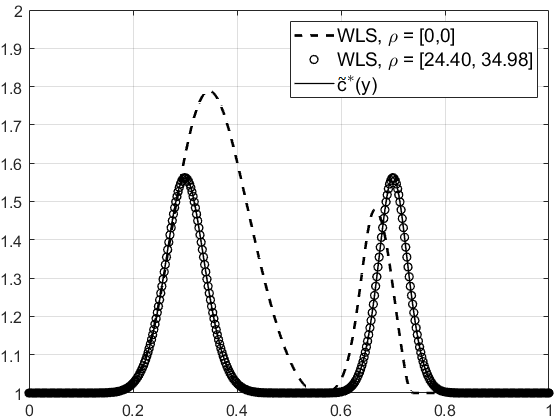}} \quad 
\subfloat[ Error function $\widetilde{c}_{comp}(y)-c^{\ast}(y)$ for the WLS with an optimal choice of $\rho = \rho^{\ast}(l)$ (solid line), compared to the Runge-Kutta method (dashed line).]{\includegraphics[width
=.45\textwidth]{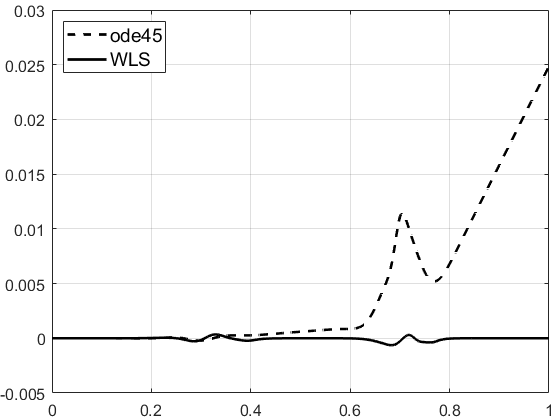}}
\end{center}
\caption{\emph{The horizontal axis depicts the values of $y$, see the text of section \ref{sec:7.2} for an explanation. (a) Computed functions and  (b) reconstruction error $\protect\widetilde{c}_{comp}(y)-\protect%
\widetilde{c}^{\ast}(y)$ for WLS and MATLAB built-in "ode45", a six-stage, fifth-order, Runge-Kutta method, where $\protect\widetilde{c}_{comp}(y)$ and $\protect\widetilde{c}^{\ast}(y)$
denote the computed and true functions respectively. See Figure \ref{fig4} (b) for the values of the optimal parameter  $\rho^{\ast}(l)$. }}
\label{fig6}
\end{figure}

Assume that the function $r_{comp}(x)$ attains negative values on the
intervals 
\begin{equation*}
(a_1,a_2),(a_3,a_4),\dots,(a_{n},a_{n+1})
\label{a}
\end{equation*}
and positive values on 
\begin{equation*}
(b_1,b_2),(b_3,b_4),\dots,(b_{m},b_{m+1}),
\label{b}
\end{equation*}%
where numbers $a_1,a_2,\dots$ and $b_1,b_2,\dots$ are found approximately.
To illustrate our idea, consider the simplest case of the function $r(x)$, depicted on Figure \ref{fig4} (a). The corresponding function $\widetilde{p}(x)$ is defined in (\ref{7.10}). We explain below how do we find the values $a_1,a_2,\dots$ and $b_1,b_2,\dots$. We first, solve the problem (\ref{7.91})-(\ref{7.92}) by the Runge-Kutta method on the subinterval $x \in (a_1,a_2)$ without using the condition $\widetilde{p}''(0) = 0$, as we have described earlier. Next, we solve problem (\ref{7.91})-(\ref{7.92}) on the interval $x \in (b_1,b_2)$ by WLS with a particular choice of the weight function. This approach provides a higher accuracy, compared with the Runge-Kutta method. The latter method results in a great deal of instability for $x \in (b_1,b_2)$, e.g. see Figure \ref{fig5} (b) and Figure \ref{fig6} (b). More precisely we minimize the following functional: 
\begin{equation}
\mathcal{L}_{\rho}(\widetilde{p}) = \int_{b_1}^{b_2} \left\vert \widetilde{p}%
^{\prime \prime }(x)\widetilde{p}(x)/2-(\widetilde{p}^{\prime}(x))^2/4 - r(x)
\right\vert^2 e^{-2 \rho (x-b_1)} dx, \hspace{0.3em} \rho \in \mathbb{R}, 
\hspace{0.3em} \rho > 1,  \label{7.11}
\end{equation}
where $\rho$ is the parameter, which is similar to the parameter $\lambda$ in the functionl $K_{\lambda,\gamma}$. Numbers $\widetilde{p}(b_1) = \widetilde{p}(a_2), \hspace{0.3em} \widetilde{p}^{\prime }(b_1) = \widetilde{p}^{\prime }(a_2), \hspace{0.3em} \widetilde{p}%
^{\prime \prime }(b_1) = \widetilde{p}^{\prime \prime }(a_2)$ are found from the solution of the problem (\ref{7.91})-(\ref{7.93}) on $x \in (a_1,a_2)$ via the Runge-Kutta method, with $a_1 = 0$ in our particular case. Simultaneously with the minimizer we find the solution to the problem. See the results obtained via the described procedure on Figures \ref{fig4}, \ref{fig6}. The initial guess for the minimization of functional (\ref{7.11}) is the cubic parabola $\widetilde{p}_0(x)$, satisfying $\widetilde{p}_0(x) = \widetilde{p}(a_2), \hspace{0.3em} \widetilde{p}_0^{\prime }(x) = \widetilde{p}^{\prime }(a_2)$, \newline $\widetilde{p}_0^{\prime \prime }(x) = \widetilde{p}^{\prime \prime }(a_2)$ at $x = b_1$ and $\widetilde{p}_0(x) = 1$ at $x = b_2$.

\vspace{0.25em}
In the more general case, assuming that $b_i = a_{i+1}$. Then the procedure is as follows:
\begin{enumerate}[label=\alph*)]
\item for each interval $(a_i,a_{i+1})$ the Runge-Kutta method is applied,
\item for each interval $(b_i,b_{i+1})$ WLS is applied.
\end{enumerate}
\vspace{0.25em}

\noindent We have established numerically that this procedure is stable as long as not too many intervals are involved.

\begin{remark}
Functional (\ref{7.11}) is very similar to the functional used in the convexification in \cite{klibanov2017globally}. In particular, a Carleman estimate for the operator $d^2/dx^2$ with the weight function $e^{-2 \rho (x-b_1)}$ can be proven. In addition, analogs of Theorems 2-6 can be proven, we refer to \cite{klibanov2017globally} for similar results.
\end{remark}

As soon as $\widetilde{p}(x)$ is computed for $x \in (0,a)$, we solve (\ref{7.93}) to find the dependence $y = y(x)$. Since the number $a$ is only an upper estimate for the interval of values of $x$ and since we know that $y \in (0,1)$, then we solve (\ref{7.93}) by the Runge-Kutta method until either $y \geq \overline{y} = 1$ or $x = a$. In our numerical experiments we observe that we always stop computation at $y = 1$. The latter allows us to tabulate the values of $\widetilde{c}(y)$ via $\widetilde{p}(x) = \widetilde{c}(y(x))$.
This is how Figures \ref{fig5} and \ref{fig6} are drawn.

We can choose the parameter $\rho$ that leads to the lowest relative error
in the reconstruction of $\widetilde{c}(y)$. This parameter depends on the length $l = b_2 - b_1$ of the interval $(b_1,b_2)$. For a given pair $(\widetilde{c}(y), l)$ an optimal value of the parameter $\rho(\widetilde{c}(y), l) = \rho^{\ast}(\widetilde{c}(y),l)$ is found numerically. The dependence of the optimal value 
on $l$ for two sets of functions $\widetilde{c}(y)$ is shown on Figure \ref{fig4} (b). This technique can be used to reconstruct an arbitrary function $\widetilde{c}(y)$, such that corresponding $\widetilde{p}(x)$ changes its sign on a finite number of intervals. In this paper, we use a simple model (\ref{7.10}) to simulate the data. The optimal value of the parameter $\rho^{\ast}(\widetilde{c}^{\ast}(y),l)$ depends on the function $\widetilde{c}^{\ast}(y)$. Since this paper considers only the functions of the form (\ref{7.10}), then the values of the optimal parameters can be approximated using the values obtained from the simulated data for $\widetilde{c}^{\ast}(y)$ with a given $A_c$ and $w^{\ast} \approx l$, see Figure \ref{fig4} (b). We have found the optimal value $\rho^{\ast}(\widetilde{c}^{\ast}(y),l)$ for the data simulated for two types of functions of the form (\ref{7.10}) with max$(\widetilde{c}(y)) = 1.56$ and max$(\widetilde{c}(y)) = 3.80$. Next, we have found numerically that if max$(\widetilde{c}(y)) \in [1.56,3.80]$ then the value of $\rho^{\ast}(\widetilde{c}(y),l)$ can be taken as
\begin{align*}
    &\rho^{\ast}(\widetilde{c}(y),l) = (\rho^{\ast}(\widetilde{c}^{\ast}_{1.56},l) + \rho^{\ast}(\widetilde{c}^{\ast}_{3.80},l))/2, \text{ see Figure \ref{fig5} (b)}, \\
    &\rho^{\ast}(\widetilde{c}^{\ast}_{1.56},l) = 2.1457/l, \quad \rho^{\ast}(\widetilde{c}^{\ast}_{3.80},l) = 2.1081/l + 14.40.
\end{align*}
Combining the above steps, we arrive at Algorithm \ref{alg 2}, which finds approximate values of $\widetilde{c}(y)$ on the whole interval $y \in (0,1)$.

\label{sec:7.2}

\begin{algorithm}
\caption{Reconstruction of $c(y)$ from the $r_{comp}(x)$ using WLS}
\begin{algorithmic}[1]
    \State\, Given the function $r_{comp}(x)$ computed via Algorithm \ref{alg 1}, determine the intervals $(a_1,a_2),(a_2,a_3), \dots, (a_n,a_{n+1})$, where $r_{comp}(x) > 0$ and $(b_1,b_2),(b_2,b_3), \dots, (b_m,b_{m+1})$, where $r_{comp}(x) \leq 0$.
    \State\, Solve (\ref{7.91})-(\ref{7.92}) via the Runge-Kutta method for $x \in (a_n,a_{n+1})$, $n = 1,2, \dots$.
    \State\, Compute the boundary conditions for $\widetilde{p}(b_{n}),\widetilde{p}'(b_{n}),\widetilde{p}''(b_n)$ for $b_n = a_{n+1}$ and an initial guess for $\widetilde{p}(x)$, via cubic extrapolation on $(b_n,b_{n+1})$.
    \State\, Minimize $ \mathcal{L}_{\rho}(\widetilde{p})$ on $(b_{m},b_{m+1})$, given $\rho = \rho^{\ast}(l)$, starting with $m = 1$.
    \State\, Compute the boundary conditions for $\widetilde{p}(a_{n+2}),\widetilde{p}'(a_{n+2}),\widetilde{p}''(a_{n+2})$ at $a_{n+2} = b_{m+1}$.
    \State\, Compute the function $y(x)$ given $\widetilde{p}(x)$ on $(b_{m},b_{m+1})$ found at step 4.
    \State\, \label{step fin} Repeat steps 2-6 for $n+1, m+1$ until $y(a_{n+1}) \geq 1$ or $y(b_{m+1}) \geq 1$.
\end{algorithmic}
\label{alg 2}
\end{algorithm}
\begin{remark}
It is clear from Figures \ref{fig4} (a) and \ref{fig6} (a) that the wrong choice of parameter $\rho(l)$ leads to significant reconstruction errors.
\end{remark}
\begin{remark}
Results depicted on Figure \ref{fig6} demonstrate a significant accuracy improvement of the numerical reconstruction of $\widetilde{c}(y)$ by WLS method as compared with the conventional Runge-Kutta "ode45" method of MATLAB. 
\end{remark}
\begin{center}
\textbf{Table 2.} Values of dielectric constants, computed from experimental data. References containing intervals with correct values of dielectric constants are listed in the first column. 
\begin{tabular}{l|l|l|l|l}
Object                  & max$(c_{comp})$ & $c_{bcgr}$ & $\varepsilon_r$ & $\varepsilon_{table}$ \\ \hline
1. Bush, see \cite{chuah1995dielectric}                       & 6.27          & 1       & 6.27       & [3,20]        \\ \hline
2. Wood stake, see  \cite{Table}               &  3.21        & 1       &   3.21     &  [2,6]       \\ \hline
3. Metal box (buried), see  \cite{Table}       & 4.12   & $[3,5]$ & $[12.36,20.60]$    & [10,30]        \\ \hline
4. Metallic cylinder (buried), see \cite{Table} & 5.39   & $[3,5]$ & $[16.17,26.95]$    &  [10,30]        \\ \hline
5. Plastic cylinder (buried), see \cite{Table} & 0.26   & $[3,5]$ & $[0.78,1.30]^{\ast}$    &  [0.5,1.5]    \\ \hline
\end{tabular}
\end{center}
$^{\ast}$the value of the dielectric constant of the Object 5 is taken as $\varepsilon_r$ = min$(c_{comp}c_{bcgr})$.
\vspace{0.5em}

\begin{remark}
Note that the computed values of are close to the ones of \cite{Kuzhuget:IP2012}, presented in Table 1 of section 8.
\end{remark}

\subsection{Numerical tests on experimental data}

\label{sec:7.3}

The experimental data were collected by A. Sullivan and L. Nguyen for five (5) objects. First two objects, bush and wood stake, were above the ground. The three buried objects were metallic box, plastic cylinder and metal cylinder. The data for each object contain 80 temporal samples with the time step $\Delta t = 0.133 ns$, which corresponds to imaged distance from $0$ to $3.15$ meters. We measure the function $g_0(t)=u(0,t)$. As to the function $g_1(t)=u_{x}(0,t)$, we calculate it using $g_0(t)$ and the absorbing boundary condition (\ref{2.21}) of Lemma 1 for $x_2 = 0$. The exact locations of targets were not of an interest to us since their horizontal coordinates were delivered quite accurately by GPS. For a more detailed description of the data acquisition scheme see section 7 of \cite{Kuzhuget:IP2012}.

The dielectric constants were not measured when the data were collected. However, it was known that the burial depth of every underground target was a few centimeters. Additionally, we knew which targets were located above the ground and which ones were buried in the ground. Thus, we have computed the values of dielectric constants and compared with those listed in other references. The background for Objects 3-5, was dry sand, where dielectric constant belongs to [3,5], see \cite{Table}. Note that tables usually contain ranges of values of dielectric constants rather than their exact values \cite{Table}. Thus, the problem of the determination of the material from the experimental data becomes a problem to find the interval for the values of dielectric constant of the studied object and a further comparison with the known range of values. 

\begin{figure}[htb!]
\begin{center}
\subfloat[The experimental data for Object 2 (wood stake, standing on the ground) with its envelope, see the text of the section 7.3 and \cite{Kuzhuget:IP2012} for the explanation of the negative values of the envelope.]{\includegraphics[width
=.45\textwidth]{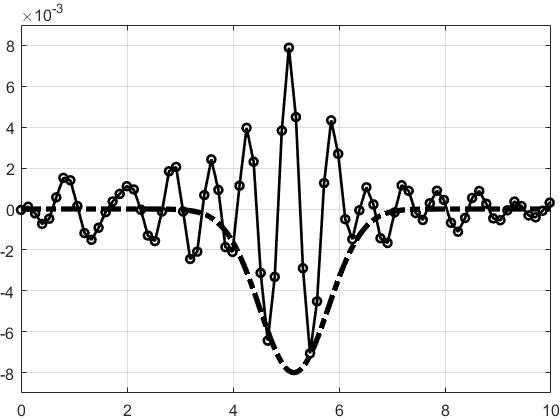}} \quad 
\subfloat[The experimental data for Object 5 (plastic cylinder, buried in sand) with its envelope, see the text of the section 7.3 and \cite{Kuzhuget:IP2012} for the explanation of the positive values of the envelope.]{\includegraphics[width =.45\textwidth]{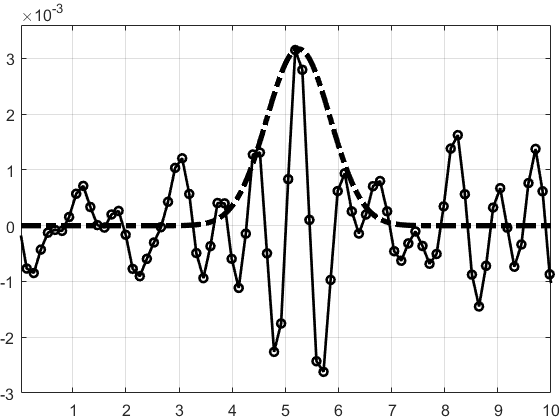}} \\
\subfloat[$\varepsilon_r(y)$, computed via WLS with optimal choice of $\rho = \rho^{\ast}(l)$ for Object 5 (plastic cylinder, buried in sand).]{\includegraphics[width
=.45\textwidth]{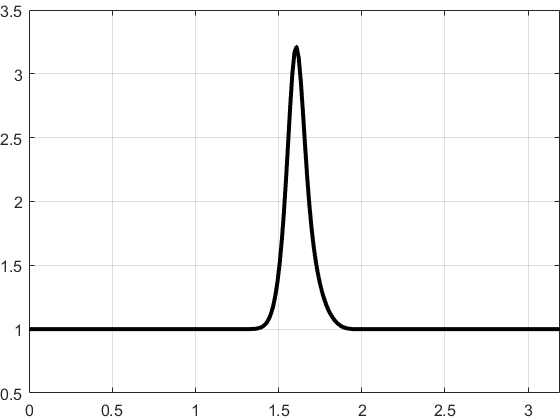}} \quad 
\subfloat[$\varepsilon_r(y)$, computed via WLS with optimal choice of $\rho = \rho^{\ast}(l)$ for Object 5 (plastic cylinder, buried in sand).]{\includegraphics[width =.45\textwidth]{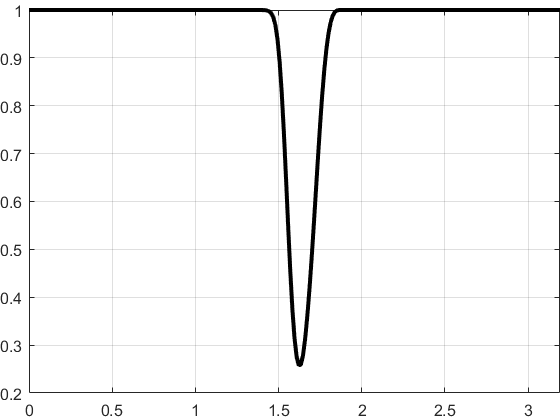}}
\end{center}
\caption{\emph{Experimental data sets for wood stake (a), standing above on the ground and plastic cylinder, buried in dry sand (b). The horizontal axis for (a)-(b) depicts time in nanoseconds, the horizontal axis for (c)-(d) depicts the distance from the antenna in meters.}}
\label{fig7}
\end{figure}

We also note that the data are 1D, whereas the objects are 3D. Thus, the images obtained here (see Figure \ref{fig7} (c),(d)) can only be considered an approximate models of the real objects. We still apply the technique of section 7.3 of \cite{Kuzhuget:IP2012} for determining of whether the positive or negative values of the data should be used for the envelope for the buried objects, see Figures \ref{fig7} (a),(b).
Figures \ref{fig7} (a),(b) depicts the data collected for two out of five objects. The data were scaled as in \cite{Kuzhuget:IP2012}, we have multiplied the data by the scaling factor $SN = 10^{-7}$.

Since our method works on the domain different from the one for which the experimental data were collected, we scale the real time $t$ in nanoseconds as $t' = 0.19\times10^9\times t$ on the interval $t' \in (0,2)$ and use the notation $t$ for the scaled time onward. For Figures \ref{fig7} (a),(b) the time was scaled back to the real time. We apply Algorithm \ref{alg 1} to reconstruct the function $r(x)$ for each object. Next, Algorithm \ref{alg 2} is applied to find the corresponding dielectric constants. The estimated values of dielectric constants are listed in the Table 2. $c_{bcgr}$ is the value of the dielectric constant of the background (air for Objects 1-2, dry sand for Objects 3-5). $\varepsilon_r$ = max$(c_{comp}c_{bcgr})$  is the estimated dielectric constant of the Object. On the other hand, if the value of the dielectric constant of the Object is less that the one of the background, e.g. Object 5 (plastic cylinder), then $\varepsilon_r$ = min$(c_{comp}c_{bcgr})$ is taken as the dielectric constant of the Object. $\varepsilon_{table}$ is the value of the dielectric constant found in the corresponding reference; those references are listed in the first column of Table 2. 

\textbf{Conclusions.} It is clear from columns 4, 5 of Table 2 that computed dielectric  constants belong to appropriate intervals, furthermore we notice that our computed dielectric constants are close to the ones of Table 1 of \cite{Kuzhuget:IP2012}.

\appendix
\section{Data simulation and pre-processing}

Let the function $c(y) \in C^3(\mathbb{R})$ be given. Consider the domain $G
= \{(y,t) \in (-\widetilde{y},\widetilde{y}) \times (0,\widetilde{T})\}, 
\hspace{0.3em} \widetilde{y} = 1.1, \widetilde{T} = 2$ and a $G^{h}$ with $%
N_y = 1600$, $N_t = 3200$. Then, to simulate the data for the inverse
problem CIP2, we solve the following initial value problem problem 
\begin{align*}
&c\left( y\right) u_{tt}=u_{yy}, \hspace{0.3em} \left( y,t\right) \in G , \\
&c(y) = 1 , \hspace{0.3em} y \in (-\widetilde{y},0] \cup [1,\widetilde{y}),
\\
\ &u\left( y,0\right) =0, \hspace{0.3em} u_{t}\left( y,0\right) = \text{exp}%
(-10^6 y^2), \\
&u_y(\widetilde{y},t) + u_t(\widetilde{y},t) = 0, \hspace{0.3em} u_y(-%
\widetilde{y},t) -u_t(-\widetilde{y},t) = 0,
\end{align*}
which is an analog of problem (\ref{2.3})-(\ref{2.4}) supplemented with the
absorbing boundary condition. In addition, for our numerical studies we
replaced the initial condition $u_{t}\left( y,0\right) =\delta(y)$ with the
gaussian $u_{t}\left( y,0\right) =$ exp$\left(-10^6 y^2\right)$. The
numerical solution of the problem stated above via finite differences
delivers the approximate function $u(y,t)$, which is used to compute
functions $g_0(t) = u(0,t)$, $g_1(t) = u_y(0,t), \hspace{0.3em} \forall t
\in (0,\widetilde{T})$.

By (\ref{3.50}), the data for the inverse problem, i.e. $s_0(t), s_1(t)$,
contain first and second derivatives of the functions $g_0(t),g_1(t)$, which
can not be retrieved numerically due to high oscillations in the values of
the considered functions. Thus, we use their envelopes to approximate the
functions $s_0(t), s_1(t)$. First, we filter the data by taking only
negative values and truncating to zero the values of the functions $%
g_0(t),g_1(t)$ which are less than $0.1\times$max$(\vert g_0(t) \vert)$ and $%
0.1\times$max$(\vert g_1(t) \vert)$ respectively. This is true for objects
having the dielectric constants greater than the one of the background. On
the other hand, for the above case of the plastic cylinder, buried in the
ground we use the positive values of the functions $g_0(t),g_1(t)$ with the
same truncation rule, see section 7.3 of \cite{Kuzhuget:IP2012}. Next, we
approximate so truncated functions with the envelopes of the following
forms: 
\begin{equation}
\widetilde{g}_{i}(t) = \pm \widetilde{g}^1_{i} exp(-\widetilde{g}^2_{i}(t-%
\widetilde{g}^3_{i})^2), \quad i = 0,1,  \label{envelope}
\end{equation}
where the sign is chosen according to the previous arguments. The parameters 
$\widetilde{g}^1_{i},\widetilde{g}^2_{i},\widetilde{g}^3_{i}$ are found using weighted least-squares curve fit ("fit" function) of the Curve Fitting Toolbox of MATLAB. Samples of approximations (\ref{envelope}) are depicted
on Figure \ref{fig2}.

Then, based on (\ref{3.50}), we set $s_0(t) = \widetilde{g}_0^{\prime }(t)$
and $s_1(t) = \widetilde{g}_0^{\prime \prime }(t)+\widetilde{g}_1^{\prime
}(t)$, where the derivatives of $\widetilde{g}_0(t), \widetilde{g}_1(t)$ are
found explicitly.

\bibliographystyle{siamplain}
\bibliography{mybib}
\end{document}